\ifx \labelsONmargin\undefined

  \ifx\RequirePackage\undefined
    \expandafter\ifx\csname amsart.sty\endcsname\relax
        \documentstyle[12pt,amssymb,amscd]{amsart}
    \fi

    \def\atSign{@@}

    \def\mathbb{\Bbb}
    \def\mathfrak{\frak}
    \def\mathbf{\bold}
    \ifx\boldsymbol\undefined	
      \def\boldsymbol#1{{\bold #1}}
    \fi

    \def\mathbit{\boldsymbol}

    \makeatletter
    \newenvironment{proof}{%
         \@ifnextchar[{%
                       \expandafter\let\expandafter\end@proof
                         \csname endpf*\endcsname
                         \my@proof
                      }{\let\end@proof\endpf\pf}%
        }{\end@proof}
    \def\my@proof[#1]{\@nameuse{pf*}{#1}}
    \makeatother

    \def\xrightarrow[#1]#2{@>{#2}>{#1}>}
    \def\xleftarrow[#1]#2{@<{#2}<{#1}<}

    \def\providecommand#1{\def#1}
    \def\emph#1{{\em #1}}
    \def\textbf#1{{\bf #1}}

    \def\mathring{\overset{\,\,{}_\circ}}
  \else
      \ifx\usepackage\RequirePackage
      \else
        \documentclass[12pt]{amsart}
        \usepackage{amssymb,amscd}
      \fi
      \ifx \mathring\undefined		
        \DeclareMathAccent{\mathring}{\mathalpha}{operators}{"17}
      \fi

      \long\def\FAKEendPROOF{\endtrivlist}
      \ifx\endproof\FAKEendPROOF
	  \def\endproof{\qed\endtrivlist}
      \fi

      \ifx\mathbit\undefined	
        \DeclareMathAlphabet{\mathbit}{OML}{cmm}{b}{it}
      \fi

      \def\atSign{@}

      \advance\oddsidemargin -0.1\textwidth
      \evensidemargin\oddsidemargin
      \def\Sb#1\endSb{_{\substack{#1}}}
      \def\Sp#1\endSp{^{\substack{#1}}}

  \fi
\makeatletter
\@ifundefined{DeclareFontShape} 
     {\@ifundefined{selectfont}
             {
                \message{We need AmS-LaTeX, this version of LaTeX does not 
                        support it}\@@end
             }{
                \def\mathcal{\cal}
                
                \def\pcyr{%
                        \def\default@family{UWCyr}%
                        \let\oldSl@\sl
                        \def\sl{\def\default@shape{it}\oldSl@}%
                        \cyracc
                        \language\Russian\family{UWCyr}\selectfont
                }
             }}{
                \DeclareFontEncoding{OT2}{}{\noaccents@}
                \DeclareFontSubstitution{OT2}{cmr}{m}{n}
                \DeclareFontFamily{OT2}{cmr}{\hyphenchar\font45 }
                \DeclareFontShape{OT2}{cmr}{m}{n}{%
                     <5><6><7><8><9><10>gen*wncyr %
                     <10.95><12><14.4><17.28><20.74><24.88> wncyr10 %
                }{}
                \DeclareFontShape{OT2}{cmr}{m}{it}{%
                     <5><6><7><8><9><10> gen * wncyi%
                     <10.95><12><14.4><17.28><20.74><24.88> wncyi10%
                }{}
                \DeclareFontShape{OT2}{cmr}{bx}{n}{%
                     <5><6><7><8><9><10> gen * wncyb%
                     <10.95><12><14.4><17.28><20.74><24.88> wncyb10%
                }{}
                \DeclareFontShape{OT2}{cmr}{m}{sl}{%
                     <-> ssub * cmr/m/it%
                }{}
                \DeclareFontShape{OT2}{cmr}{m}{sc}{%
                     <5><6><7><8><9><10>%
                     <10.95><12><14.4><17.28><20.74><24.88> wncysc10%
                }{}
                \DeclareFontFamily{OT2}{cmss}{\hyphenchar\font45 }
                \DeclareFontShape{OT2}{cmss}{m}{n}{%
                     <8><9><10> gen * wncyss%
                     <10.95><12><14.4><17.28><20.74><24.88> wncyss10%
                }{}
                
                \def\cyrencodingdefault{OT2}
                \def\pcyr{%
                        \cyracc
                        \let\encodingdefault\cyrencodingdefault
                        \language\Russian\fontencoding{OT2}\selectfont
                }
             }   

\@ifundefined{theorembodyfont} 
     {
        \def\theorembodyfont#1{\relax}
        \makeatletter
          \let\@@th@plain\th@plain
          \def\th@plain{ \@@th@plain \slshape }
        \makeatother
        \let\normalshape\relax
     }{}   

\ifx\cprime\undefined
     \def\cprime{$'$}
\fi
\makeatletter
  \def\@sect@my#1#2#3#4#5#6[#7]#8{%
\ifnum #2>\c@secnumdepth
   \let\@svsec\@empty
 \else
   \refstepcounter{#1}%
\edef\@svsec{\ifnum#2<\@m
             \@ifundefined{#1name}{}{\csname #1name\endcsname\ }\fi
\noexpand\rom{\csname the#1\endcsname.}\enspace}\fi
 \@tempskipa #5\relax
 \ifdim \@tempskipa>\z@ 
   \begingroup #6\relax
   \@hangfrom{\hskip #3\relax\@svsec}{\interlinepenalty\@M #8\par}%
   \endgroup
   \if@article\else\csname #1mark\endcsname{%
        \ifnum \c@secnumdepth >#2\relax\csname the#1\endcsname. \fi#7}\fi
\ifnum#2>\@m \else
       \let\@tempf\\ \def\\{\protect\\}\addcontentsline{toc}{#1}%
{\ifnum #2>\c@secnumdepth \else
             \protect\numberline{%
               \ifnum#2<\@m
               \@ifundefined{#1name}{}{\csname #1name\endcsname\ }\fi
               \csname the#1\endcsname.}\fi
           #8}\let\\\@tempf
     \fi
 \else
  \def\@svsechd{#6\hskip #3\@svsec
    \@ifnotempty{#8}{\ignorespaces#8\unskip
       \ifnum\spacefactor<1001.\fi}%
        \ifnum#2>\@m \else
          \let\@tempf\\ \def\\{\protect\\}\addcontentsline{toc}{#1}%
            {\ifnum #2>\c@secnumdepth \else
              \protect\numberline{%
                \ifnum#2<\@m
                \@ifundefined{#1name}{}{\csname #1name\endcsname\ }\fi
                \csname the#1\endcsname.}\fi
             #8}\let\\\@tempf\fi}%
 \fi
\@xsect{#5}}
  \ifx\RequirePackage\undefined
  \let\@sect\@sect@my             
  \fi
  \ifx\RequirePackage\undefined	
  \def\th@remark@my{\theorempreskipamount6\p@\@plus6\p@
    \theorempostskipamount\theorempreskipamount
    \def\theorem@headerfont{\it}\normalshape}
    \let\th@remark\th@remark@my
  \else				
    \let\o@@remark\th@remark

    \ifx\theorempostskipamount\undefined		
    \else						
      \def\th@remark{\o@@remark
	\ifdim\theorempostskipamount < 2pt\relax
	  \theorempostskipamount\theorempreskipamount
	     \multiply\theorempostskipamount\tw@
	     \divide\theorempostskipamount\thr@@
	\fi
      }
    \fi
  \fi
\let\myLabel\@gobble
\def\labelsONmargin{\@mparswitchfalse\def\myLabel##1{\@bsphack\marginpar
                                  {\normalshape\tiny\rm Label ##1}\@esphack}}
\ifx \url\undefined
  \def\url#1{{\tt #1}}%
\fi
\def\cyracc{\def\u##1{
                \if \i##1\char"1A%
                \else \if I##1\char"12%
                \else \accent"24 ##1\fi\fi }%
\def\"##1{\if e##1{\char"1B}%
                \else \if E##1{\char"13}%
                \else \accent"7F ##1\fi\fi }%
\def\9##1{\if##1z\char"19 
\else\if##1Z\char"11 
\else\if##1E\char"03 
\else\if##1e\char"0B 
\else\if##1u\char"18 
\else\if##1U\char"10 
\else\if##1A\char"17 
\else\if##1a\char"1F 
\else\if##1p\char"7E 
\else\if##1P\char"5E 
\else\if##1Q\char"5F 
\else\if##1q\char"7F 
\else\if##1i\char"1A 
\else\if##1I\char"12 
\else\if##1N\char"7D 
\fi
\fi
\fi
\fi
\fi
\fi
\fi
\fi
\fi
\fi
\fi
\fi
\fi
\fi
\fi
}%
\def\cydot{{\kern0pt}}}%
\def\cydot{$\cdot$}

\ifx\Russian\undefined
        \def\Russian{0\relax
    \message{Don't know the hyphenation rules for Russian^^J
                        Please do INITeX with `input  russhyph' in the 
                        command line}%
                \gdef\Russian{0\relax}%
        }
\fi
\ifx\RequirePackage\undefined			
  \def\@putname#1#2#3#4{\def\@@ref{#3}\let\old@bf\bf
        \def\bf##1{\old@bf\if?\noexpand##1?{#4}\else##1\fi}%
	#1{#2}%
        \let\bf\old@bf}
\else						
  \def\@putname#1#2#3#4{\def\@@ref{#3}\let\old@bf\bf	
	\let\old@reset@font\reset@font			
        \def\bf##1{\old@bf\if?\noexpand##1?{#4}\else##1\fi}%
	\def\reset@font##1##2{\old@reset@font##1\if?\noexpand##2?{#4}\else##2\fi}#1{#2}%
        \let\bf\old@bf\let\reset@font\old@reset@font}
\fi
\let\my@ref=\ref
\def\ref#1{\@putname\my@ref{#1}{#1}{\tiny\rm\@@ref}}
\let\my@pageref=\pageref
\def\pageref#1{\@putname\my@pageref{#1}{#1}{\tiny\rm\@@ref}}
\let\my@cite=\cite
\def\cite#1{\@putname\my@cite{#1}{\@citeb}{\tiny\rm\@@ref}}
\makeatother
\fi
\relax 
\theoremstyle{plain} 
\theorembodyfont{\sl}

\numberwithin{equation}{section}
\theorembodyfont{\rm}
\theoremstyle{definition}

\newtheorem{definition}{Definition}[section]
\newtheorem{conjecture}[definition]{Conjecture}
\newtheorem{example}[definition]{Example}

\theoremstyle{remark}

\newtheorem{remark}[definition]{Remark} 

\theoremstyle{plain} 
\theorembodyfont{\sl}

\newtheorem{theorem}[definition]{Theorem}
\newtheorem{lemma}[definition]{Lemma}
\newtheorem{corollary}[definition]{Corollary}

\begin{document}
\bibliographystyle{amsplain}
\relax 

\title{ On associated variety for Lie superalgebras }

\author{ Michel Duflo and Vera Serganova }

\date{ \today }

\address{ Dept. of Mathematics, University of California at Berkeley,
Berkeley, CA 94720 }

\email{serganov\atSign{}math.berkeley.edu}

\maketitle
\begin{abstract}

We define the associated variety $ X_{M} $ of a module $ M $ over a
finite-dimensional superalgebra $ {\mathfrak g} $, and show how to extract information
about $ M $ from these geometric data. $ X_{M} $ is a subvariety of the cone $ X $ of
self-commuting odd elements.

For finite-dimensional $ M $, $ X_{M} $ is invariant under the action of the
underlying Lie group $ G_{0} $. For simple superalgebra with invariant symmetric
form, $ X $ has finitely many $ G_{0} $-orbits; we associate a number (rank) to each
such orbit. One can also associate a number (degree of atypicality) to an
irreducible finite-dimensional representation.

We prove that if $ M $ is an irreducible $ {\mathfrak g} $-module of degree of atypicality
$ k $, then $ X_{M} $ lies in the closure of all orbits on $ X $ of rank $ k $. If $ {\mathfrak g}={\mathfrak g}{\mathfrak l}\left(m|n\right) $
we prove that $ X_{M} $ coincides with this closure.

\end{abstract}
\tableofcontents

\section{Introduction }

In this paper we introduce a notion of associated variety for a module
over a Lie superalgebra. This is a superanalogue of an associated variety of
Harish-Chandra modules. Associated varieties have many interesting
applications in classical representation theory (see, for example,
\cite{Vog91Ass,KatOchi01Deg,NishOchiTan01Ber}).

The associated variety for a Lie superalgebra is a subvariety of a cone
$ X\subset{\mathfrak g}_{1} $ of self-commuting odd elements. This cone $ X $ was studied by Caroline
Gruson, see \cite{Gru00Ide,Gru00Coh,Gru03Coh}. She used geometric properties
of $ X $ to obtain important results about cohomology of Lie superalgebras.

While the associated variety in classical representation theory is
trivial if a module is finite-dimensional, finite-dimensional modules over
classical Lie superalgebras have interesting associated varieties. Since
finite-dimensional representation theory of superalgebras still has many
open problems, we hope that the associated variety will have some
application in this theory. In particular, it should help to describe
analytic properties of supercharacters and cohomolgy groups. Some simple
applications are given in Sections 3 and 7.

Let us outline the results of this paper. In Section 2 we give a
definition and formulate simple properties of associated variety. In Section
3 we construct a coherent sheaf on $ X $ associated with $ M $ and prove a
criterion of projectivity for modules over certain
Lie superalgebras. In Section 4 we discuss geometry of $ X $. Section 5 contains
main theorems (Theorem~\ref{th2} and Theorem~\ref{th3}) about the associated
varieties for simple classical
contragredient superalgebras. In Section 6 we prove Theorem~\ref{th2}.
In Section 7 we give some applications of Theorem~\ref{th2} to
supercharacters. Sections 8,~9 and 10 contain a proof of Theorem~\ref{th3}.

\section{Definition and basic properties }

Let $ {\mathfrak g}={\mathfrak g}_{0}\oplus{\mathfrak g}_{1} $ be a finite-dimensional complex Lie superalgebra, $ G_{0} $
denote a simply-connected connected Lie group with Lie algebra $ {\mathfrak g}_{0} $. Let
\begin{equation}
X=\left\{x\in{\mathfrak g}_{1} \mid \left[x,x\right]=0\right\}.
\notag\end{equation}
It is clear that $ X $ is $ G_{0} $-invariant Zariski closed cone in $ {\mathfrak g}_{1} $. Let $ M $ be a
$ {\mathfrak g} $-module. For each $ x\in X $ put $ M_{x}=\operatorname{Ker} x/xM $ and define
\begin{equation}
X_{M}=\left\{x\in X \mid M_{x}\not=0\right\}.
\notag\end{equation}
We call $ X_{M} $ the {\em associated variety\/} of $ M $.

\begin{lemma} \label{lm1}\myLabel{lm1}\relax  If $ M $ is a finite-dimensional $ {\mathfrak g} $-module, then $ X_{M} $ is Zarisky
closed $ G_{0} $-invariant subvariety.

\end{lemma}

\begin{proof} Since $ M $ is finite-dimensional, $ M $ is a $ G_{0} $-module. For each
$ g\in G_{0} $ and $ x\in M $ one has
\begin{equation}
M_{\operatorname{Ad}_{g}\left(x\right)}=gM_{x},
\notag\end{equation}
that implies Lemma.\end{proof}

\begin{lemma} \label{lm2}\myLabel{lm2}\relax 
\begin{enumerate}
\item
If $ M=U\left({\mathfrak g}\right)\otimes_{U\left({\mathfrak g}_{0}\right)}M_{0} $ for some $ {\mathfrak g}_{0} $-module $ M_{0} $, then $ X_{M}=\left\{0\right\} $;
\item
If $ M={\mathbb C} $ is trivial, then $ X_{M}=X $;
\item
For any $ {\mathfrak g} $-modules $ M $ and $ N $, one has $ X_{M\oplus N}=X_{M}\cup X_{N} $;
\item
For any $ {\mathfrak g} $-modules $ M $ and $ N $, one has $ X_{M\otimes N}=X_{M}\cap X_{N} $;
\item
For any finite-dimensional $ {\mathfrak g} $-module $ M $, $ X_{M^{*}}=X_{M} $;
\item
For any finite-dimensional $ {\mathfrak g} $-module $ M $ and any $ x\in X $, $ \operatorname{sdim} M=\operatorname{sdim} M_{x} $.
\end{enumerate}
\end{lemma}

\begin{proof} Properties 2,3,5 follow directly from definition. To prove 1,
let $ x\in X $ and $ x\not=0 $. Let $ \left\{v_{j}\right\}_{j\in J} $ be a basis of $ M_{0} $ and $ x_{1},\dots ,x_{m} $ be a basis of $ {\mathfrak g}_{1} $
such that $ x=x_{1} $. Then by PBW for Lie superalgebras $ x_{i_{1}}x_{i_{2}}\dots x_{i_{k}}\otimes v_{j} $ for all
$ 1\leq i_{1}<i_{2}<\dots <i_{k}\leq m $, $ j\in J $ form a basis of $ M $. The action of $ x=x_{1} $ in this basis is
easy to write and it is clear that $ \operatorname{Ker} x=xM $ is spanned by the vectors
$ x_{1}x_{i_{2}}\dots x_{i_{k}}\otimes v_{j} $.

Now let us show (4). We will prove that $ M_{x}=0 $ implies $ \left(M\otimes N\right)_{x}=0 $.
Indeed, $ M_{x}=0 $ implies that $ M $ is a free $ {\mathbb C}\left[x\right] $-module. Tensor product of a
free $ {\mathbb C}\left[x\right] $-module with any $ {\mathbb C}\left[x\right] $-module is free. Therefore $ M\otimes N $ is free over
$ {\mathbb C}\left[x\right] $ and $ \left(M\otimes N\right)_{x}=0 $.

Finally we will prove (6). Let $ \Pi\left(N\right) $ stand for the superspace
isomorphic to $ N $ with switched parity. Since $ M/\operatorname{Ker} x $ is isomorphic to $ \Pi\left(xM\right) $,
then
\begin{equation}
\operatorname{sdim} M=\operatorname{sdim} \operatorname{Ker} x+\operatorname{sdim} \Pi\left(x M\right)=\operatorname{sdim} \operatorname{Ker} x-\operatorname{sdim} x M=\operatorname{sdim}\left(\operatorname{Ker} x /xM\right)=\operatorname{sdim}
M_{x}.
\notag\end{equation}
\end{proof}

\section{Localization and projective modules }

Let $ {\mathcal O}_{X} $ denote the structure sheaf of $ X $. Then $ {\mathcal O}_{X}\otimes M $ is the sheaf of
sections of a trivial vector bundle with fiber isomorphic to $ M $. Let $ \partial:{\mathcal O}_{X}\otimes M \to
{\mathcal O}_{X}\otimes M $ be the map defined by
\begin{equation}
\partial\varphi\left(x\right)=x\varphi\left(x\right)
\notag\end{equation}
for any $ x\in X $, $ \varphi\in{\mathcal O}_{X}\otimes M $. Clearly $ \partial^{2}=0 $ and the cohomology $ {\mathcal M} $ of $ \partial $ is a
quasi-coherent sheaf on $ X $. If $ M $ is finite-dimensional, then $ {\mathcal M} $ is coherent.

For any $ x\in X $ denote by $ {\mathcal O}_{x} $ the local ring of $ x $, by $ {\mathcal I}_{x} $ the maximal ideal.
Then the fiber $ {\mathcal M}_{x} $ is the the cohomology of $ \partial:{\mathcal O}_{x}\otimes M \to {\mathcal O}_{x}\otimes M $. The evaluation
map $ j_{x}:{\mathcal O}_{x}\otimes M \to M $ satisfies $ j_{x}\circ\partial=x\circ j_{x} $. Hence we have the maps
\begin{equation}
j_{x}:\operatorname{Ker} \partial \to\operatorname{Ker} x\text{, }j_{x}:\operatorname{Im} \partial \to xM.
\notag\end{equation}
One can easily check that the latter map is surjective. Therefore $ j_{x} $ induces
the map $ \bar{j}_{x}\colon {\mathcal M}_{x} \to M_{x} $, and $ \operatorname{Im} \bar{j}_{x}\cong{\mathcal M}_{x}/{\mathcal I}_{x}{\mathcal M}_{x} $.

\begin{lemma} \label{lm3}\myLabel{lm3}\relax  Let $ M $ be a finite-dimensional $ {\mathfrak g} $-module. The support of $ {\mathcal M} $ is
contained in
$ X_{M} $. The map $ \bar{j}_{x} $ is surjective for a generic point $ x\in X $. In particular,
if $ X_{M}=X $, then $ \operatorname{supp} {\mathcal M}=X $.
\end{lemma}

\begin{proof} First, we will show that for any $ x\in X\backslash X_{M} $ there exists a
neighborhood $ U $ of $ x $ such that $ {\mathcal M}\left(U\right)=0 $. Indeed, there exists a map
$ i_{x}:M \to M $ such that $ x\circ i_{x}=id. $ Therefore in some neighborhood $ U $ of $ x $ there
exists a map $ i:{\mathcal O}\left(U\right)\otimes M \to {\mathcal O}\left(U\right)\otimes M $ such that $ \partial\circ i=\operatorname{id} $ and $ i\left(x\right)=i_{x} $, hence $ {\mathcal M}\left(U\right)=0 $.
Thus, $ x $ does not belong to the support of $ {\mathcal M} $ and we have obtained that $ \operatorname{supp}
{\mathcal M}\subset X_{M} $.

To prove the second statement let $ x\in X $ be such that $ \dim  xM $ is maximal
possible. Let $ m\in\operatorname{Ker} x $. Then there exists some neighborhood $ U $ of $ x $ and
$ \varphi\in{\mathcal O}\left(U\right)\otimes M $ such that $ \partial\varphi=0 $ and $ \varphi\left(x\right)=m $. By definition $ \varphi\in{\mathcal M}_{x} $ and $ \bar{j}_{x}\left(\varphi\right)=m $.\end{proof}

\begin{corollary} \label{cor59}\myLabel{cor59}\relax  Let $ x\in X $ be a generic point, then in some neighborhood $ U $
of $ x $, the sheaf $ {\mathcal M}_{U} $ coincides with the sheaf of section of a vector bundle
with fiber $ M_{x} $.

\end{corollary}

Let $ X_{M}\not=X $. Then $ {\mathcal M} $ is the extension by zero of the sheaf $ {\mathcal M}_{X_{M}} $. If we denote
by $ {\mathcal M}\left(x\right) $ the image of $ \bar{j}_{x} $, then $ {\mathcal M}_{X_{M}} $ locally is the sheaf of sections of the
vector bundle with fiber $ {\mathcal M}\left(x\right) $ for a generic $ x\in X_{M} $. Note that
$ {\mathcal M}\left(x\right)\subset M_{x} $, but usually $ {\mathcal M}\left(x\right)\not=M_{x} $, as one can see from the following example.

\begin{example} \label{ex1}\myLabel{ex1}\relax  Let $ {\mathfrak g}={\mathfrak s}{\mathfrak l}\left(1|n\right) $. Then $ {\mathfrak g}_{1}={\mathfrak g}\left(-1\right)\oplus{\mathfrak g}\left(1\right) $, where $ {\mathfrak g}\left(-1\right) $ and $ {\mathfrak g}\left(1\right) $ are
abelian superalgebras. Assume that $ M $ is a typical irreducible $ {\mathfrak g} $-module.
Then $ X_{M}=\left\{0\right\} $, $ M_{0}=M $ and $ {\mathcal M}\left(0\right)=M^{{\mathfrak g}\left(1\right)}\oplus M^{{\mathfrak g}\left(-1\right)} $.

\end{example}

Let $ {\mathcal F} $ be the category of finite-dimensional $ {\mathfrak g} $-modules semisimple over
$ {\mathfrak g}_{0} $. The latter condition is automatic if $ {\mathfrak g}_{0} $ is semisimple.

\begin{theorem} \label{th142}\myLabel{th142}\relax  Assume that $ {\mathfrak g}_{0} $ is a reductive Lie algebra and elements
of $ X $ span $ {\mathfrak g}_{1} $. Then $ M\in{\mathcal F} $ is projective iff $ X_{M}=\left\{0\right\} $.

\end{theorem}

\begin{proof} Let $ M $ be projective. Since $ M $ is a quotient of $ U\left({\mathfrak g}\right)\otimes_{U\left({\mathfrak g}_{0}\right)}M $,
then $ M $ is a direct summand of $ U\left({\mathfrak g}\right)\otimes_{U\left({\mathfrak g}_{0}\right)}M $. By Lemma~\ref{lm2} $ \left(1\right) $ and $ \left(3\right) X_{M}=\left\{0\right\} $.

To prove the assertion in opposite direction we need the following
lemma. Let $ H_{\text{red}}^{\cdot}\left({\mathfrak g},M\right) $ denote the cohomology of $ {\mathfrak g} $, induced by cocycles
trivial on the center of $ {\mathfrak g}_{0} $.

\begin{lemma} \label{lm142}\myLabel{lm142}\relax  Let $ {\mathfrak g} $ satisfy the condition of Theorem, $ M\in{\mathcal F} $ and $ X_{M}=\left\{0\right\} $.
Then $ H_{\text{red}}^{1}\left({\mathfrak g},M\right)=\left\{0\right\} $.

\end{lemma}

\begin{proof} Let $ \varphi\in{\mathfrak g}^{*}\otimes M $ be a $ 1 $-cocycle. Then for any $ x\in X $ we have $ x\varphi\left(x\right)=0 $.
Thus, $ \varphi $ induces a global section of $ {\mathcal M} $. Since $ \varphi\left(0\right)=0 $ and the $ \operatorname{supp} {\mathcal M}=\left\{0\right\} $ by
Lemma~\ref{lm3}, this global section must be zero. Therefore there exists $ \psi\left(x\right) $
such that $ x\psi\left(x\right)=\varphi\left(x\right) $ for all $ x\in X $. But $ \varphi $ is a linear function, therefore $ \psi $ is
constant. If $ d $ is the differential in the cohomology complex, $ \eta=\varphi-d\psi $ is a
$ 1 $-cocycle homologically equivalent to $ \varphi $. On the other hand, $ \eta\left(x\right)=0 $ for any
$ x\in X $, and since elements of $ X $ span $ {\mathfrak g}_{1} $, we have $ \eta\left({\mathfrak g}_{1}\right)=0 $. The restriction $ \eta $ on
$ {\mathfrak g}_{0} $ is a $ 1 $-cocycle for a Lie algebra $ {\mathfrak g}_{0} $. But $ {\mathfrak g}_{0} $ is reductive, hence
$ H_{\text{red}}^{1}\left({\mathfrak g}_{0},M\right)=0 $ and therefore $ \eta=d\nu $. We have shown that $ \varphi $ induces the trivial
cohomology class. Thus, $ H_{\text{red}}^{1}\left({\mathfrak g},M\right)=\left\{0\right\} $.\end{proof}

Now, assume that $ X_{M}=\left\{0\right\} $. We have to show that $ \operatorname{Ext}^{1}\left(M,N\right)=\left\{0\right\} $, the latter
is equivalent to $ H_{\text{red}}^{1}\left({\mathfrak g}, M^{*}\otimes N\right)=\left\{0\right\} $. By Lemma~\ref{lm2} $ \left(4\right),\left(5\right) $ we have $ X_{M^{*}\otimes N}=\left\{0\right\} $.
Therefore $ H_{\text{red}}^{1}\left({\mathfrak g}, M^{*}\otimes N\right)=\left\{0\right\} $ and $ M $ is projective.\end{proof}

\begin{remark} \label{rem79}\myLabel{rem79}\relax  Note that the conditions of Theorem~\ref{th142} hold for any
simple classical superalgebra except $ {\mathfrak o}{\mathfrak s}{\mathfrak p}\left(1|2n\right) $. In case
$ {\mathfrak g}={\mathfrak o}{\mathfrak s}{\mathfrak p}\left(1|2n\right) $, $ X=\left\{0\right\} $ and $ {\mathcal F} $ is semi-simple, hence every finite-dimensional
module is projective. In general, however, Theorem~\ref{th142} is not true if we
drop the assumption that $ X $ spans $ {\mathfrak g}_{1} $. Indeed, let $ {\mathfrak g}={\mathfrak q}\left(1\right) $, in other words $ {\mathfrak g} $
has a basis of an even element $ C $ and an odd element $ T $ such that $ \left[T,T\right]=C $.
Then $ X=\left\{0\right\} $, but not every module in $ {\mathcal F} $ is projective. For example, the
trivial one-dimensional $ {\mathfrak g} $-module is not projective.

\end{remark}

\section{The structure of $ X $ for contragredient simple Lie superalgebras }

Let $ {\mathfrak g} $ be a contragredient finite-dimensional Lie superalgebra with
indecomposable Cartan matrix, i.e. $ {\mathfrak g} $ is
isomorphic to one from the following list: $ {\mathfrak s}{\mathfrak l}\left(m|n\right) $ if $ m\not=n $, $ {\mathfrak g}{\mathfrak l}\left(n|n\right) $,
$ {\mathfrak o}{\mathfrak s}{\mathfrak p}\left(m|2n\right) $, $ D\left(\alpha\right) $, $ F_{4} $ or $ G_{3} $ (for definitions see \cite{Kac77Lie}).

\begin{remark} \label{rem1}\myLabel{rem1}\relax  The Lie superalgebras we consider are simple except one case. For a
simple Lie superalgebra $ {\mathfrak p}{\mathfrak s}{\mathfrak l}\left(n|n\right) $ the Cartan matrix
is degenerate and we consider the corresponding Kac-Moody Lie
superalgebra which is isomorphic to $ {\mathfrak g}{\mathfrak l}\left(n|n\right) $. Later we will do the proofs
for $ {\mathfrak g}{\mathfrak l}\left(m|n\right) $ even if $ m\not=n $, in this case $ {\mathfrak g}{\mathfrak l}\left(m|n\right)\cong{\mathfrak s}{\mathfrak l}\left(m|n\right)\oplus{\mathbb C} $.

\end{remark}

We fix a Cartan subalgebra $ {\mathfrak h}\subset{\mathfrak g} $. In this case the Cartan subalgebra of $ {\mathfrak g} $
coincides with a Cartan subalgebra of $ {\mathfrak g}_{0} $ and $ {\mathfrak g} $ has a root decomposition
\begin{equation}
{\mathfrak g}={\mathfrak h}\oplus \oplus_{\alpha\in\Delta}{\mathfrak g}_{\alpha},
\notag\end{equation}
each root space $ {\mathfrak g}_{\alpha} $ is one dimensional. The parity of $ \alpha\in\Delta $ by definition is
equal to the parity of the root space $ {\mathfrak g}_{\alpha} $. The invariant bilinear form $ \left(\cdot,\cdot\right) $
on $ {\mathfrak h}^{*} $ is not positive definite and some of odd roots are isotropic. For a
non-isotropic $ \beta $ we denote by $ \beta^{\vee} $ the element of $ {\mathfrak h} $ such that
$ \alpha\left(\beta^{\vee}\right)=\frac{2\left(\alpha,\beta\right)}{\left(\beta,\beta\right)} $. Let $ S $ denote the set of subsets of mutually orthogonal
linearly independent isotropic roots of $ \Delta_{1} $, i.e. an element of $ S $ is
$ {\text A}=\left\{\alpha_{1},\dots ,\alpha_{k} \mid \left(\alpha_{i},\alpha_{j}\right)=0\right\} $. The Weyl group $ W $ of $ {\mathfrak g}_{0} $ acts on $ S $ in the obvious
way. Put $ S_{k}=\left\{{\text A}\in S \mid |{\text A}|=k\right\} $, here $ S_{0}=\left\{\varnothing\right\} $.

\begin{theorem} \label{th1}\myLabel{th1}\relax  There are finitely many $ G_{0} $-orbits on $ X $. These orbits
are in one-to one correspondence with $ W $-orbits in $ S $.

\end{theorem}

\begin{proof} We define the map $ \Phi:S \to X/G_{0} $ in the following way. Let
$ {\text A}=\left\{\alpha_{1},\dots ,\alpha_{k}\right\}\in S $, choose a non-zero $ x_{i}\in{\mathfrak g}_{\alpha_{i}} $ and put $ x=x_{1}+\dots +x_{k}\in X $. By
definition $ \Phi\left({\text A}\right)=G_{0}x $. To see that $ \Phi\left({\text A}\right) $ does not depend on a choice of $ x_{i} $
note that since $ \alpha_{1},\dots ,\alpha_{k} $ are linearly independent, for any other choice
\begin{equation}
x'=\Sigma x_{i}'=\Sigma c_{i}x_{i}
\notag\end{equation}
there is $ h\in{\mathfrak h} $ such that $ c_{i}=e^{\alpha_{i}\left(h\right)} $ and therefore
\begin{equation}
x'=\exp \left(\operatorname{ad}\left(h\right)\right)\left(x\right).
\notag\end{equation}
If $ {\text B}=w\left({\text A}\right) $ for some $ w\in W $, then clearly $ \Phi\left({\text B}\right) $ and $ \Phi\left({\text A}\right) $ belong to the same
orbit. Therefore $ \Phi $ induces the map $ \bar{\Phi}:S/W \to X/G_{0} $. We check
case by case that $ \bar{\Phi} $ is injective and surjective.

If $ {\mathfrak g} $ is $ {\mathfrak s}{\mathfrak l}\left(m|n\right) $ or $ {\mathfrak g}{\mathfrak l}\left(n|n\right) $, $ {\mathfrak g} $ has a natural $ {\mathbb Z} $ grading
$ {\mathfrak g}={\mathfrak g}\left(-1\right)\oplus{\mathfrak g}\left(0\right)\oplus{\mathfrak g}\left(1\right) $ such that $ {\mathfrak g}_{0}={\mathfrak g}\left(0\right) $, $ {\mathfrak g}_{1}={\mathfrak g}\left(1\right)\oplus{\mathfrak g}\left(-1\right) $. The orbits of $ W $ on $ S $ are
enumerated by the pairs of numbers $ \left(p,q\right) $, where $ p=|{\text A}\cap\Delta\left({\mathfrak g}\left(1\right)\right)| $,
$ q=|{\text A}\cap\Delta\left({\mathfrak g}\left(-1\right)\right)| $. The orbits of $ G_{0} $ on $ X $ are enumerated by the same pairs of
numbers $ \left(p,q\right) $ in the following way. If $ x=x^{+}+x^{-} $, where $ x^{\pm}\in{\mathfrak g}\left(\pm1\right) $, then
$ p=\operatorname{rank}\left(x^{+}\right) $, $ q=\operatorname{rank}\left(x^{-}\right) $. We can see by the construction of $ \bar{\Phi} $, that $ \bar{\Phi} $ maps
$ \left(p,q\right) $-orbit on $ S $ to the $ \left(p,q\right) $-orbit on $ X $.

Let $ {\mathfrak g}={\mathfrak o}{\mathfrak s}{\mathfrak p}\left(m|2n\right) $. If $ m=2l+1 $ or $ m=2l $ with $ l>n $, then the $ W $-orbits on $ S $ are
in one-to-one correspondence with $ \left\{0,1,2,\dots ,\min \left(l,n\right)\right\} $. Namely, $ {\text A} $ and $ {\text B} $ are
on the same orbit if they have the same number of elements. As it was shown
in \cite{Gru00Ide}, $ X $ can be identified with the set of all linear maps $ x:{\mathbb C}^{m} \to
{\mathbb C}^{2n} $, such that $ \operatorname{Im} x $ is an isotropic subspace in $ {\mathbb C}^{2n} $ and $ \operatorname{Im} x^{*} $ is an
isotropic subspace in $ {\mathbb C}^{m} $. Furthermore, $ x,y\in X $ belong to the same $ G_{0} $-orbit iff
$ \operatorname{rank}\left(x\right)=\operatorname{rank}\left(y\right) $. One can see that rank $ \Phi\left({\text A}\right)=|{\text A}| $.

Now let $ {\mathfrak g}={\mathfrak o}{\mathfrak s}{\mathfrak p}\left(2l|2n\right) $ where $ l\leq n $. If $ {\text A},{\text B}\in S $ and $ |{\text A}|=|{\text B}|<l $, then $ {\text A} $ and $ {\text B} $
are on the same $ W $-orbit. In the same way if $ \operatorname{rank}\left(x\right)=\operatorname{rank}\left(y\right)<l $, then $ x $ and $ y $
are on the same $ G_{0} $-orbit. However, the set of all $ x\in{\mathfrak g}_{1} $ of maximal rank
splits in two orbits, since the Grassmannian of maximal isotropic subspaces
in $ {\mathbb C}^{2l} $ has two connected components. In the same way $ S_{l} $ splits in two
$ W $-orbits. Hence in this case again $ \bar{\Phi} $ is a bijection.

If $ {\mathfrak g} $ is one of exceptional Lie superalgebras $ D\left(\alpha\right) $, $ G_{3} $ or $ F_{4} $, then
the direct calculation shows that $ X $ has two $ G_{0} $-orbits: $ \left\{0\right\} $ and the orbit
of a highest vector in $ {\mathfrak g}_{1} $. The set $ S $ also consists of two $ W $-orbits: $ \varnothing $ and
the set of all isotropic roots in $ \Delta $.\end{proof}

\begin{remark} \label{rem393}\myLabel{rem393}\relax  Note that in our situation the representation of $ G_{0} $ in $ {\mathfrak g}_{1} $
is symplectic and multiplicity free (see \cite{Knop05inv}). The cone $ X $ is the
preimage of 0 under the moment map $ {\mathfrak g}_{1} \to {\mathfrak g}_{0}^{*} $.

\end{remark}

We use the notation $ \Phi:S \to X/G_{0} $ introduced in the proof of Theorem~%
\ref{th1}. Using the explicit description of $ G_{0} $-orbits on $ X $ and the description of
roots systems, which can be found in \cite{Kac77Lie}, one can check
the following statements case by case. We omit this checking here.

\begin{lemma} \label{lm110}\myLabel{lm110}\relax  Let $ {\text A},{\text B}\in S. $
\begin{enumerate}
\item
If $ \alpha\in\Delta $ is a linear combination of roots from $ {\text A} $, then $ \alpha\in{\text A}\cup-{\text A} $;
\item
If $ |{\text A}|\leq|{\text B}| $, then there exists $ w\in W $ such that $ w\left({\text A}\right)\subset{\text B}\cup-{\text B} $;
\item
$ \Phi\left({\text A}\right) $ lies in the closure of $ \Phi\left({\text B}\right) $ iff $ w\left({\text A}\right)\subset{\text B} $ for some $ w\in W $.
\end{enumerate}
\end{lemma}

By $ {\text A}^{\perp} $ we denote the set of all weights orthogonal to $ {\text A} $ with respect to
the standard form on $ {\mathfrak h}^{*} $.

\begin{theorem} \label{th7}\myLabel{th7}\relax  Let $ {\text A}\in S $. Then $ \dim  \Phi\left({\text A}\right)=\frac{|\Delta_{1}\backslash{\text A}^{\perp}|}{2}+|{\text A}| $.

\end{theorem}

\begin{proof} Let $ {\text A}=\left\{\alpha_{1},\dots ,\alpha_{k}\right\} $, $ x=x_{1}+\dots +x_{k} $ for some choice of $ x_{i}\in{\mathfrak g}_{\alpha_{i}} $,
$ y=y_{1}+\dots +y_{k} $ for some $ y_{i}\in{\mathfrak g}_{-\alpha_{i}} $. Let $ h=\left[x,y\right] $, $ h_{i}=\left[x_{i},y_{i}\right] $. Clearly,
$ h=h_{1}+\dots +h_{k} $ and $ h,x,y $ generate the $ {\mathfrak s}{\mathfrak l}\left(1|1\right) $-subalgebra in $ {\mathfrak g} $. With respect
to this subalgebra $ {\mathfrak g} $ has a decomposition
\begin{equation}
{\mathfrak g}=\oplus_{\mu}{\mathfrak g}^{\mu},
\notag\end{equation}
where
\begin{equation}
{\mathfrak g}^{\mu}=\left\{g\in{\mathfrak g} \mid \left[h,g\right]=\mu g\right\}.
\notag\end{equation}
Note that
\begin{equation}
\dim  \left[{\mathfrak g},x\right]=\sum_{\mu}\dim  \left[{\mathfrak g}^{\mu},x\right],
\notag\end{equation}
and from the description of irreducible $ {\mathfrak s}{\mathfrak l}\left(1|1\right) $-modules for $ \mu\not=0 $
\begin{equation}
\dim  \left[{\mathfrak g}^{\mu},x\right]=\frac{\dim  {\mathfrak g}^{\mu}}{2}.
\notag\end{equation}
On the other hand, for $ \mu\not=0 \operatorname{sdim} {\mathfrak g}^{\mu}=0 $, and therefore
\begin{equation}
\dim  {\mathfrak g}^{\mu}=2 \dim  {\mathfrak g}_{1}^{\mu}.
\notag\end{equation}
Observe that for a generic choice of $ x_{i}\in{\mathfrak g}_{\alpha_{i}} $, $ {\mathfrak g}_{\beta}\subset{\mathfrak g}^{0} $ iff $ \left(\beta,\alpha_{i}\right)=0 $ for
all $ i\leq k $. Indeed, for generic
choice of $ x_{i} $ the condition $ \beta\left(h\right)=0 $ implies $ \beta\left(h_{i}\right)=0 $ for all $ i $, and
therefore $ \left(\beta,\alpha_{i}\right)=0 $ for all $ i $. Hence
\begin{equation}
\oplus_{\mu\not=0}{\mathfrak g}_{1}^{\mu}=\oplus_{\alpha\in\Delta_{1}\backslash{\text A}^{\perp}}{\mathfrak g}_{\alpha}
\notag\end{equation}
and
\begin{equation}
\sum_{\mu\not=0} \dim  \left[{\mathfrak g}^{\mu},x\right]=\sum_{\mu\not=0} \dim  {\mathfrak g}_{1}^{\mu}=|\Delta_{1}\backslash{\text A}^{\perp}|.
\notag\end{equation}
To calculate $ \dim  \left[{\mathfrak g}^{0},x\right] $ note that
\begin{equation}
{\mathfrak g}^{0}={\mathfrak h}\oplus\oplus_{\beta\in\Delta\cap{\text A}^{\perp}}{\mathfrak g}_{\beta}.
\notag\end{equation}
We claim that
\begin{equation}
\left[{\mathfrak g}^{0},x\right]=\oplus_{i=1}^{k}{\mathbb C}h_{i}\oplus\oplus_{i=1}^{k}{\mathfrak g}_{\alpha_{i}},
\notag\end{equation}
hence $ \dim  \left[{\mathfrak g}^{0},x\right]=2k $. Indeed, if $ \left(\beta,\alpha_{i}\right)=0,\beta\not=\pm\alpha_{i} $ then $ \beta\pm\alpha_{i}\notin\Delta $. Therefore
$ \left[x,{\mathfrak g}_{\beta}\right]=0 $ for any $ \beta\in\Delta\cap{\text A}^{\perp},\beta\not=-\alpha_{i} $. Furthermore, $ \left[x,{\mathfrak g}_{-\alpha_{i}}\right]={\mathbb C}h_{i} $ and $ \left[x,{\mathfrak h}\right]=\oplus_{i=1}^{k}{\mathfrak g}_{\alpha_{i}} $.
Thus, we obtain
\begin{equation}
\dim  \left[{\mathfrak g},x\right]=|\Delta_{1}\backslash{\text A}^{\perp}|+2k.
\label{equ99}\end{equation}\myLabel{equ99,}\relax 
Now the statement will follow from the lemma.

\begin{lemma} \label{lm71}\myLabel{lm71}\relax  $ \operatorname{sdim} \left[{\mathfrak g},x\right]=0 $.

\end{lemma}

\begin{proof} Define the odd skew-symmetric form on $ {\mathfrak g} $ by
\begin{equation}
\omega\left(y,z\right)=\left(x,\left[y,z\right]\right).
\notag\end{equation}
Obviously the kernel of $ \omega $ coincides with centralizer $ C_{{\mathfrak g}}\left(x\right) $. Thus, $ \omega $ is
non-degenerate odd skew-symmetric form on $ C_{{\mathfrak g}}\left(x\right) $. Hence $ \operatorname{sdim} {\mathfrak g}/C_{{\mathfrak g}}\left(x\right)=0 $. But
$ \left[{\mathfrak g},x\right]\cong\Pi\left({\mathfrak g}/C_{{\mathfrak g}}\left(x\right)\right) $, which implies the lemma.\end{proof}

Lemma implies that $ \dim \left[{\mathfrak g}_{0},x\right]=1/2 \dim \left[{\mathfrak g},x\right] $. Since $ \dim  G_{0}x=\dim \left[{\mathfrak g}_{0},x\right] $, the
theorem follows from~\eqref{equ99}.

{}\end{proof}

\begin{corollary} \label{cor73}\myLabel{cor73}\relax  If $ |{\text A}|=|{\text B}| $, then $ \dim  \Phi\left({\text A}\right)= \dim  \Phi\left({\text B}\right) $.

\end{corollary}

\begin{proof} Follows from Theorem~\ref{th1} and Lemma~\ref{lm110} (2).\end{proof}

The maximal number of isotropic mutually orthogonal linearly
independent roots is called the {\em defect\/} of $ {\mathfrak g} $. This notion was introduced in
\cite{KacWak94Aff}. One can see that the
defect of $ {\mathfrak g} $ is equal to the dimension of maximal isotropic subspace in $ {\mathfrak h}^{*} $.
All exceptional Lie superalgebras has defect 1. The defect of $ {\mathfrak s}{\mathfrak l}\left(m|n\right) $ is
$ \min \left(m,n\right) $, the defect $ {\mathfrak o}{\mathfrak s}{\mathfrak p}\left(2l+1|2n\right) $ and $ {\mathfrak o}{\mathfrak s}{\mathfrak p}\left(2l|2n\right) $ is $ \min  \left(l,n\right) $.

\begin{corollary} \label{cor71}\myLabel{cor71}\relax  Let $ d $ be the defect of $ {\mathfrak g} $. Then the irreducible
components of $ X $ are in bijection with $ W $-orbits on $ S_{d} $. If all odd roots of $ {\mathfrak g} $
are isotropic, then the dimension of each component equals $ \frac{\dim 
{\mathfrak g}_{1}}{2}=\frac{|\Delta_{1}|}{2} $.

\end{corollary}

\begin{proof} As follows from Theorem~\ref{th1} and Lemma~\ref{lm110} (3), each
irreducible component is the closure of $ \Phi\left({\text A}\right) $ for a maximal $ {\text A}\in S $. By Lemma~%
\ref{lm110} $ \left(2\right) |{\text A}|=d $. Hence the first statement. Theorem~\ref{th7} immediately
implies the statement about dimension.\end{proof}

\begin{corollary} \label{cor72}\myLabel{cor72}\relax  If all odd roots of $ {\mathfrak g} $ are isotropic, then the
codimension of $ \Phi\left({\text A}\right) $ in $ X $ equals $ \frac{|\Delta_{1}\cap{\text A}^{\perp}|}{2}-|{\text A}| $.

\end{corollary}

\begin{proof} The codimension equals $ \dim  X-\dim  \Phi\left({\text A}\right) $. Using Theorem~\ref{th7} and
Corollary~\ref{cor71}
\begin{equation}
\operatorname{codim} \Phi\left({\text A}\right)=\frac{|\Delta_{1}|-|\Delta_{1}\backslash{\text A}^{\perp}|}{2}-|{\text A}|=\frac{|\Delta_{1}\cap{\text A}^{\perp}|}{2}-|{\text A}|.
\notag\end{equation}
\end{proof}

\section{Central character and the main theorems }

Let us fix a Borel subalgebra $ {\mathfrak b}\subset{\mathfrak g} $ by choosing a decomposition
$ \Delta=\Delta^{+}\cup\Delta^{-} $. Note that this choice is not unique but our consideration will
not depend on it. Later we will use different Borel subalgebras in some
proofs. Let
\begin{equation}
\rho=\frac{1}{2}\sum_{\alpha\in\Delta_{0}^{+}}\alpha-\frac{1}{2}\sum_{\alpha\in\Delta_{1}^{+}}\alpha,
\notag\end{equation}
and define the shifted action of $ W $ on $ {\mathfrak h}^{*} $ by
\begin{equation}
\lambda^{w}=w\left(\lambda+\rho\right)-\rho.
\notag\end{equation}
By $ M_{\lambda} $ we denote the Verma module $ U\left({\mathfrak g}\right)\otimes_{U\left({\mathfrak b}\right)}C_{\lambda} $, and by $ L_{\lambda} $ we denote the
unique irreducible quotient of $ M_{\lambda} $. We say that $ \lambda\in{\mathfrak h}^{*} $ is {\em integral dominant\/}
if $ L_{\lambda} $ is finite-dimensional. We denote by $ \Sigma $ the set of all integral
dominant weights.

Let $ Z $ denote the center of the universal enveloping algebra $ U\left({\mathfrak g}\right) $. One
can see that any $ z\in Z $ acts as a scalar $ \chi_{\lambda}\left(z\right) $ on $ M_{\lambda} $ and $ L_{\lambda} $. Therefore $ \lambda\in{\mathfrak h}^{*} $
defines a central character $ \chi_{\lambda}:Z \to {\mathbb C} $. Let
\begin{equation}
{\mathfrak h}_{\chi}=\left\{\mu\in{\mathfrak h}^{*} \mid \chi_{\mu}=\chi\right\}.
\notag\end{equation}

\begin{lemma} \label{lm6}\myLabel{lm6}\relax  Let $ \chi=\chi_{\lambda} $, $ {\text A}\in S $ be a maximal set of linearly independent
mutually orthogonal isotropic roots orthogonal to $ \lambda+\rho $ and $ {\mathfrak t}_{\lambda}=\lambda+\oplus_{\alpha\in{\text A}}{\mathbb C}\alpha $. Then
\begin{equation}
{\mathfrak h}_{\chi}=\bigcup_{w\in W} {\mathfrak t}_{\lambda}^{w}.
\notag\end{equation}
\end{lemma}

\begin{proof} Easily follows from the description of the $ Z $ formulated in
\cite{Kac84Lap} and proven \cite{Ser99Inv} and in \cite{Gor04Kac}.\end{proof}

Let us fix a central character $ \chi $. For each $ \lambda\in{\mathfrak h}_{\chi} $ define $ S_{\lambda}\subset S $ by the
following
\begin{equation}
S_{\lambda}=\left\{{\text A}\in S \mid {\text A}\subset\left(\lambda+\rho\right)^{\perp}\right\}.
\notag\end{equation}
Put
\begin{equation}
S_{\chi}=\cup_{\lambda\in{\mathfrak h}_{\chi}}S_{\lambda}.
\notag\end{equation}

\begin{lemma} \label{lm5}\myLabel{lm5}\relax  There exists a number $ k $ such that $ S_{\chi}=\bigcup_{i\leq k}S_{i} $.

\end{lemma}

\begin{proof} It follows easily from Lemma~\ref{lm6} that $ S_{\chi} $ is $ W $-invariant.
Furthermore, if $ {\text A}\in S_{\chi} $ and $ {\text A}' $ is obtained from $ {\text A} $ by multiplication of some
roots in $ {\text A} $ on $ -1 $, then $ {\text A}'\in S_{\chi} $. Hence the statement follows from Lemma~\ref{lm110}
(1) and (2).\end{proof}

The number $ k $ is called the {\em degree of atypicality\/} of
$ \chi $. The degree of atypicality of $ \lambda $ is by definition the degree of atypicality
of $ \chi_{\lambda} $. If $ k=0 $, then $ \chi $ is called {\em typical}. It is clear
that the degree of atypicality of $ \chi $ is not bigger than the defect of $ {\mathfrak g} $.

Let $ X_{k}=\Phi\left(S_{k}\right) $, $ \bar{X}_{k} $ denote the closure of $ X_{k} $. Lemma~\ref{lm110} $ \left(3\right) $ implies
that
\begin{equation}
\bar{X}_{k}=\bigcup_{i=0}^{k}\Phi\left(X_{i}\right).
\notag\end{equation}

\begin{theorem} \label{th2}\myLabel{th2}\relax  Let $ {\mathfrak g} $ be a contragredient simple Lie superalgebra, $ M $ be a
$ {\mathfrak g} $-module which admits central character $ \chi $, the degree of atypicality of $ \chi $
be equal to $ k $. Then $ X_{M}\subset\bar{X}_{k} $.

\end{theorem}

\begin{theorem} \label{th3}\myLabel{th3}\relax  Let $ {\mathfrak g}={\mathfrak g}{\mathfrak l}\left(m|n\right) $ or $ {\mathfrak s}{\mathfrak l}\left(m|n\right) $. For any
integral dominant $ \lambda\in{\mathfrak h}^{*} $ with degree of atypicality $ k $, $ X_{L_{\lambda}}=\bar{X}_{k} $.

\end{theorem}

\begin{conjecture} \label{con1}\myLabel{con1}\relax  Let $ {\mathfrak g} $ be a contragredient simple Lie superalgebra.
For any integral dominant $ \lambda\in{\mathfrak h}^{*} $ with degree of atypicality $ k $, $ X_{L_{\lambda}}=\bar{X}_{k} $.

\end{conjecture}

First, observe that the conjecture is true for the typical character.

\begin{theorem} \label{th10}\myLabel{th10}\relax  If $ \lambda $ is typical, then $ X_{L_{\lambda}}=\left\{0\right\} $.

\end{theorem}

\begin{proof} If $ \lambda $ is typical, then $ L_{\lambda} $ is a direct summand of some induced
module $ U\left({\mathfrak g}\right)\otimes_{U\left({\mathfrak g}_{0}\right)}M_{0} $(see \cite{Kac78Rep} ). Therefore Theorem follows from Lemma~%
\ref{lm2} $ \left(1\right) $ and (3).\end{proof}

\section{The structure of a generic fiber and the proof of Theorem\protect ~\protect \ref{th2}}

In this section we discuss properties of the fiber $ M_{x} $ over a point
$ x\in X_{M} $. Let $ C_{{\mathfrak g}}\left(x\right) $ be the centralizer of $ x\in X $, then by definition
$ {\mathfrak g}_{x}=C_{{\mathfrak g}}\left(x\right)/\left[x,{\mathfrak g}\right] $.

\begin{lemma} \label{lm201}\myLabel{lm201}\relax  The subspace $ \left[x,{\mathfrak g}\right] $ is an ideal in $ C_{{\mathfrak g}}\left(x\right) $. Let $ {\mathfrak m}^{\perp} $ denote the
orthogonal complement to $ {\mathfrak m} $ with respect to the invariant form on $ {\mathfrak g} $. Then
$ \left[x,{\mathfrak g}\right]^{\perp}=C_{{\mathfrak g}}\left(x\right) $.
\end{lemma}

\begin{proof} Let $ u\in C_{{\mathfrak g}}\left(x\right) $, $ v\in\left[x,{\mathfrak g}\right] $. Then $ v=\left[x,z\right] $ and
\begin{equation}
\left[u,\left[x,z\right]\right]=\left(-1\right)^{p\left(u\right)}\left[x,\left[u,z\right]\right]\in\left[x,{\mathfrak g}\right].
\notag\end{equation}
The second statement follows from the identity
\begin{equation}
\left(u,\left[x,z\right]\right)=-\left(\left[u,x\right],z\right).
\notag\end{equation}
\end{proof}

\begin{lemma} \label{lm204}\myLabel{lm204}\relax  $ M_{x} $ is a $ C_{{\mathfrak g}}\left(x\right) $-module trivial over $ \left[x,{\mathfrak g}\right] $.

\end{lemma}

\begin{proof} Let $ m\in\operatorname{Ker} x $, $ v=\left[x,z\right]\in\left[x,{\mathfrak g}\right] $. Then
\begin{equation}
v m=x z m-\left(-1\right)^{p\left(z\right)}z x m=x z m\in x M.
\notag\end{equation}
\end{proof}

In the case of contragredient finite-dimensional superalgebra we can
describe $ {\mathfrak g}_{x} $ precisely. Let $ {\text A}=\left\{\alpha_{1},\dots ,\alpha_{k}\right\}\in S $, $ x\in X $, and $ x=x_{1}+\dots +x_{k} $, where
$ x_{i}\in{\mathfrak g}_{\alpha_{i}} $, $ {\mathfrak h}_{\alpha}=\left[{\mathfrak g}_{\alpha},{\mathfrak g}_{-\alpha}\right] $. Define $ {\text A}'={\text A}^{\perp}\cap\Delta\backslash\left({\text A}\cup-{\text A}\right) $, $ {\mathfrak h}_{{\text A}}={\mathfrak h}_{\alpha_{1}}\oplus\dots \oplus{\mathfrak h}_{\alpha_{k}} $.
\begin{lemma} \label{lm202}\myLabel{lm202}\relax  If $ {\mathfrak g} $ is finite-dimensional contragredient superalgebra,
$ {\text A}=\left\{\alpha_{1},\dots ,\alpha_{k}\right\}\in S $, $ x\in X $, and $ x=x_{1}+\dots +x_{k} $, where $ x_{i}\in{\mathfrak g}_{\alpha_{i}} $. Then $ C_{{\mathfrak g}}\left(x\right) $ can be
decomposed in a semidirect sum $ {\mathfrak g}_{x}+\left[x,{\mathfrak g}\right] $, where $ {\mathfrak g}_{x} $ is spanned by the root
spaces $ {\mathfrak g}_{\alpha} $ for all $ \alpha\in{\text A}' $ and $ {\mathfrak h}_{x}\subset{\mathfrak h}_{{\text A}}^{\perp} $ is such that $ {\mathfrak h}_{x}\oplus{\mathfrak h}_{{\text A}}={\mathfrak h}_{{\text A}}^{\perp} $. Furthermore, def
$ {\mathfrak g}_{x}= $def $ {\mathfrak g}-k $.

\end{lemma}

\begin{proof} We use the same argument as in the proof of Theorem~\ref{th7}.
Let $ h $ and $ {\mathfrak g}^{\mu} $ be as in this proof. First, there is an isomorphism
\begin{equation}
{\mathfrak g}_{x}\cong{\mathfrak g}^{0}\cap C_{{\mathfrak g}}\left(x\right)/{\mathfrak g}^{0}\cap\left[x,{\mathfrak g}\right].
\notag\end{equation}
Then we notice that
\begin{equation}
{\mathfrak g}^{0}\cap C_{{\mathfrak g}}\left(x\right)={\mathfrak h}_{{\text A}}^{\perp}\oplus \oplus_{\alpha\in{\text A}^{\perp}\cap\Delta\backslash-{\text A}} {\mathfrak g}_{\alpha}\text{, }{\mathfrak g}^{0}\cap\left[x,{\mathfrak g}\right]={\mathfrak h}_{{\text A}}\oplus{\mathfrak g}_{\alpha_{1}}\oplus\dots \oplus{\mathfrak g}_{\alpha_{k}}.
\notag\end{equation}
Choose $ {\mathfrak h}_{x} $ in such a way that $ {\mathfrak g}_{x}=\left({\mathfrak h}_{x}\oplus\oplus_{\alpha\in{\text A}'} {\mathfrak g}_{\alpha} \right) $ is a subalgebra, then
\begin{equation}
{\mathfrak g}^{0}\cap C_{{\mathfrak g}}\left(x\right)={\mathfrak g}_{x}\oplus{\mathfrak g}^{0}\cap\left[x,{\mathfrak g}\right].
\notag\end{equation}
\end{proof}

\begin{remark} \label{rem203}\myLabel{rem203}\relax  If $ {\mathfrak g}={\mathfrak g}{\mathfrak l}\left(m|n\right) $, then $ {\mathfrak g}_{x}\cong{\mathfrak g}{\mathfrak l}\left(m-k|n-k\right) $. If $ {\mathfrak g}={\mathfrak o}{\mathfrak s}{\mathfrak p}\left(m|2n\right) $, then
$ {\mathfrak g}_{x}\cong{\mathfrak o}{\mathfrak s}{\mathfrak p}\left(m-2k|2n-2k\right) $. If $ {\mathfrak g}=D\left(\alpha\right) $ and $ x\not=0 $, then $ {\mathfrak g}_{x}\cong{\mathbb C} $. For $ {\mathfrak g}=G_{3} $ or $ F_{4} $ for a
non-zero $ x\in X $, $ {\mathfrak g}_{x} $ is isomorphic to $ {\mathfrak s}{\mathfrak l}\left(2\right) $ and $ {\mathfrak s}{\mathfrak l}\left(3\right) $ respectively.

\end{remark}

\begin{lemma} \label{lm205}\myLabel{lm205}\relax  Let $ x:V \to V $ be an odd linear operator such that $ x^{2}=0 $.
Assume that $ V=W\oplus U $, where $ W $ is a trivial $ {\mathbb C}\left[x\right] $-submodule and $ U $ is a free
$ {\mathbb C}\left[x\right] $-module. Let $ S\left(V\right)^{x} $ denote the space of $ x $-invariants in $ S\left(V\right) $. Then
$ S\left(V\right)^{x}=S\left(W\right)\otimes S\left(U\right)^{x} $ and $ S\left(U\right)^{x}\subset S\left(U\right)U^{x} $.

\end{lemma}

Let $ U\left({\mathfrak g}\right)^{\operatorname{ad}\left(x\right)} $ denote the subalgebra of $ \operatorname{ad}\left(x\right) $-invariants in $ U\left({\mathfrak g}\right) $, $ I_{x} $ be
the left ideal in $ U\left({\mathfrak g}\right) $ generated by $ \left[x,{\mathfrak g}\right] $. One has the following sequence
\begin{equation}
U\left({\mathfrak g}_{x}\right) \xrightarrow[]{\iota} U\left({\mathfrak g}\right)^{\operatorname{ad}\left(x\right)} \xrightarrow[]{\pi} U\left({\mathfrak g}\right)^{\operatorname{ad}\left(x\right)}/I_{x}\cap U\left({\mathfrak g}\right)^{\operatorname{ad}\left(x\right)}.
\notag\end{equation}
Let $ \phi=\pi\circ\iota $.

\begin{lemma} \label{lm206}\myLabel{lm206}\relax  The map $ \phi:U\left({\mathfrak g}_{x}\right) \to U\left({\mathfrak g}\right)^{\operatorname{ad}\left(x\right)}/I_{x}\cap U\left({\mathfrak g}\right)^{\operatorname{ad}\left(x\right)} $ is an
isomorphism of vector spaces.

\end{lemma}

\begin{proof} Since $ I_{x}\cap U\left({\mathfrak g}_{x}\right)=\left\{0\right\} $, $ \phi $ is injective. To prove surjectivity of
$ \phi $ use PBW and the corresponding sequence for symmetric algebras
\begin{equation}
S\left({\mathfrak g}_{x}\right) \to S\left({\mathfrak g}\right)^{\operatorname{ad}\left(x\right)} \to S\left({\mathfrak g}\right)^{\operatorname{ad}\left(x\right)}/J_{x}\cap S\left({\mathfrak g}\right)^{\operatorname{ad}\left(x\right)},
\notag\end{equation}
where $ J_{x}=\left[x,{\mathfrak g}\right]S\left({\mathfrak g}\right) $. Apply Lemma~\ref{lm205} with $ V={\mathfrak g} $. Then
$ W={\mathfrak g}_{x} $, $ U^{x}=\left[x,{\mathfrak g}\right] $, and we obtain $ S\left({\mathfrak g}\right)^{\operatorname{ad}\left(x\right)}=S\left({\mathfrak g}_{x}\right)\otimes S\left(U\right)^{x} $ and $ S\left(U\right)^{x}\subset\left[x,{\mathfrak g}\right]S\left(U\right) $.
Thus, $ gr \phi $ is an isomorphism. Hence $ \phi $ is an isomorphism.\end{proof}

Define the map $ \eta:U\left({\mathfrak g}\right)^{\operatorname{ad}\left(x\right)} \to U\left({\mathfrak g}_{x}\right) $ by putting $ \eta=\phi^{-1}\circ\pi $. As follows
from Lemma~\ref{lm204} for any $ u\in U\left({\mathfrak g}\right)^{\operatorname{ad}\left(x\right)} $, $ m\in M_{x} $
\begin{equation}
u m=\eta\left(u\right) m
\label{equ98}\end{equation}\myLabel{equ98,}\relax 
Note that $ \iota,\pi $ are homomorphisms of $ {\mathfrak g}_{x} $-modules (with respect to the
adjoint action). The center $ Z $ of $ U\left({\mathfrak g}\right) $ obviously is a subalgebra in
$ U\left({\mathfrak g}\right)^{\operatorname{ad}\left(x\right)} $. Let $ Z\left({\mathfrak g}_{x}\right) $ be the center of $ U\left({\mathfrak g}_{x}\right) $. Since $ \eta $ is a homomorphism of
$ {\mathfrak g}_{x} $-modules, $ \eta\left(Z\right)\subset Z\left({\mathfrak g}_{x}\right) $. We are going to describe the dual map
\begin{equation}
\eta^{*}\colon \operatorname{Hom} \left(Z\left({\mathfrak g}_{x}\right),{\mathbb C}\right) \to \operatorname{Hom} \left(Z,{\mathbb C}\right).
\notag\end{equation}
Choose a borel subalgebra $ {\mathfrak b}\subset{\mathfrak g} $ such that $ \alpha_{1},\dots ,\alpha_{k} $ are simple roots.

\begin{lemma} \label{lm207}\myLabel{lm207}\relax  Let $ \lambda\in{\mathfrak h}^{*} $ satisfy $ \left(\lambda+\rho,\alpha_{1}\right)=\dots =\left(\lambda+\rho,\alpha_{k}\right)=0 $. Then $ \left(L_{\lambda}\right)_{x}\not=0 $.
In particular the highest vector $ v $ belongs to $ \left(L_{\lambda}\right)_{x} $.

\end{lemma}

\begin{proof} Clearly, $ v\in\operatorname{Ker} x $. If $ v=x w $, then one can choose $ w $ with
weight $ \lambda-\alpha_{i} $ for some $ i $. However, $ L_{\lambda} $ does not have vectors of such weight.\end{proof}

\begin{corollary} \label{cor207}\myLabel{cor207}\relax  Let $ \lambda $ be as in Lemma~\ref{lm207} and $ \mu $ be the restriction
of $ \lambda $ to $ {\mathfrak h}_{x} $. Let $ \chi_{\mu}\in\operatorname{Hom} \left(Z\left({\mathfrak g}_{x}\right),{\mathbb C}\right) $ be induced by $ \mu $ and $ \chi_{\lambda}\in\operatorname{Hom} \left(Z,{\mathbb C}\right) $ be
induced by $ \lambda $ via Harish-Chandra homomorphism. Then $ \eta^{*}\left(\chi_{\mu}\right)=\chi_{\lambda} $.

\end{corollary}

\begin{corollary} \label{cor208}\myLabel{cor208}\relax  Let $ \chi\in\operatorname{Hom} \left(Z\left({\mathfrak g}_{x}\right),{\mathbb C}\right) $ and have the degree of
atypicality $ s $. Then the degree of atypicality of $ \eta^{*}\left(\chi\right) $ equals $ s+k $.

\end{corollary}

Corollary~\ref{cor208} implies Theorem~\ref{th2}. It also implies the following

\begin{theorem} \label{th209}\myLabel{th209}\relax  Let $ M $ admit a central character with degree of
atypicality $ k $, and $ x\in X_{k} $. Then $ {\mathfrak g}_{x} $-module $ M_{x} $ admits a typical central
character. In particular, if $ M_{x} $ is finite dimensional, it is semi-simple
over $ {\mathfrak g}_{x} $, and therefore over $ C_{{\mathfrak g}}\left(x\right) $.

\end{theorem}

\begin{theorem} \label{th210}\myLabel{th210}\relax  If $ {\mathfrak g}\not={\mathfrak o}{\mathfrak s}{\mathfrak p}\left(2l|2n\right) $ or $ D\left(\alpha\right) $, then $ \eta^{*} $ is injective, and
therefore $ \eta $ is surjective. If $ {\mathfrak g}={\mathfrak o}{\mathfrak s}{\mathfrak p}\left(2l|2n\right) $ or $ D\left(\alpha\right) $, then a preimage of $ \eta^{*} $
has at most two elements.

\end{theorem}

\begin{proof} Let $ {\text A}=\left\{\alpha_{1},\dots ,\alpha_{k}\right\} $, $ x=x_{1}+\dots +x_{k} $, $ x_{i}\in{\mathfrak g}_{\alpha_{i}} $, $ {\mathfrak b} $ be such that
$ \alpha_{1},..,\alpha_{k} $ are simple. Let
\begin{equation}
W'=\left\{w\in W \mid w\left({\text A}\right)\subset{\text A}\cup-{\text A}\right\},
\notag\end{equation}
and $ W\left({\mathfrak g}_{x}\right) $ denote the Weyl group of $ {\mathfrak g}_{x} $. Clearly, $ W\left({\mathfrak g}_{x}\right)\subset W $. One can show that if
$ {\mathfrak g}\not={\mathfrak o}{\mathfrak s}{\mathfrak p}\left(2l|2n\right) $ or $ D\left(\alpha\right) $, then $ W'=W\left({\mathfrak g}_{x}\right)\times W'' $, where $ W'' $ consists of all elements
which act trivially on $ {\mathfrak g}_{x} $.

Let $ \lambda\in{\mathfrak h}^{*} $, $ \left(\lambda+\rho,\alpha_{i}\right)=0 $ for all $ i=1,\dots ,k $. Then\footnote{We also use the fact that $ \rho_{x}=\frac{1}{2}\displaystyle\Sigma_{\alpha\in\Delta_{0}\left({\mathfrak g}_{x}\right)}\alpha-\frac{1}{2}\displaystyle\Sigma_{\alpha\in\Delta_{1}\left({\mathfrak g}_{x}\right)}\alpha=\rho_{|{\mathfrak h}_{x}} $. Hence
the shifted action of $ W\left({\mathfrak g}_{x}\right) $ is the same.}
\begin{equation}
{\mathfrak h}_{\chi_{\lambda}}\cap\left({\mathfrak h}_{{\text A}}^{\perp}\right)^{*}=\bigcup_{w\in W'} {\mathfrak t}_{\lambda}^{w}\text{, }{\mathfrak h}_{\chi_{\lambda}}\cap{\mathfrak h}_{x}=\bigcup_{w\in W'} \left({\mathfrak t}_{\lambda}^{w}\cap{\mathfrak h}_{x}^{*}\right).
\notag\end{equation}
Let $ {\mathfrak g}\not={\mathfrak o}{\mathfrak s}{\mathfrak p}\left(2l|2n\right) $ or $ D\left(\alpha\right) $ and $ \mu $ be the restriction of $ \lambda $ on $ {\mathfrak h}_{x} $. Then
\begin{equation}
{\mathfrak h}_{\chi_{\lambda}}\cap{\mathfrak h}_{x}=\bigcup_{w\in W\left({\mathfrak g}_{x}\right)} \left({\mathfrak t}_{\lambda}\cap{\mathfrak h}_{x}^{*}\right)^{w}=\bigcup_{w\in W\left({\mathfrak g}_{x}\right)} {\mathfrak t}_{\mu}^{w}=\left({\mathfrak h}_{x}\right)_{\chi_{\mu}},
\notag\end{equation}
that shows $ \left(\eta^{*}\right)^{-1}\left(\chi_{\lambda}\right)=\chi_{\mu} $.

In case $ {\mathfrak g}={\mathfrak o}{\mathfrak s}{\mathfrak p}\left(2l|2n\right) $ or $ D\left(\alpha\right) $, $ W\left({\mathfrak g}_{x}\right)\times W'' $ has index
2 in $ W' $. Take $ u\in W' $, $ u\notin W\left({\mathfrak g}_{x}\right)\times W'' $, let $ \mu $ be the restriction of $ \lambda $ on $ {\mathfrak h}_{x} $ and $ \mu' $
be the restriction of $ \lambda^{u} $ on $ {\mathfrak h}_{x} $. Then
\begin{equation}
{\mathfrak h}_{\chi_{\lambda}}\cap{\mathfrak h}_{x}=\bigcup_{w\in W'} \left({\mathfrak t}_{\lambda}\cap{\mathfrak h}_{x}^{*}\right)^{w}=\bigcup_{w\in W\left({\mathfrak g}_{x}\right)} \left({\mathfrak t}_{\mu}^{w}\cup{\mathfrak t}_{\mu'}^{w}\right)=\left({\mathfrak h}_{x}\right)_{\chi_{\mu}}\cup\left({\mathfrak h}_{x}\right)_{\chi_{\mu'}}.
\notag\end{equation}
Therefore $ \left(\eta^{*}\right)^{-1}\left(\chi_{\lambda}\right)=\left\{\chi_{\mu},\chi_{\mu'}\right\} $.\end{proof}

Assume that $ M $ is finite-dimensional and has central character $ \chi $ with
degree of atypicality $ k $. Let $ x\in\bar{X}_{k} $.
Let
\begin{equation}
Y_{x}=\left\{y\in\left({\mathfrak g}_{x}\right)_{1} \mid \left[y,y\right]=0\right\}.
\notag\end{equation}
Then
\begin{equation}
x+Y_{x}\subset X
\label{equ83}\end{equation}\myLabel{equ83,}\relax 
Define the coherent sheaf $ {\mathcal N} $ on $ Y_{x} $ as the cohomology of
\begin{equation}
\partial:{\mathcal O}_{Y_{x}}\otimes M_{x} \to {\mathcal O}_{Y_{x}}\otimes M_{x}.
\notag\end{equation}
Let $ {\mathcal N}\left(x\right) $ be the image of the fiber $ {\mathcal N}_{x} $ in $ M_{x} $ under the evaluation map.

\begin{theorem} \label{th211}\myLabel{th211}\relax  $ {\mathcal M}\left(x\right)={\mathcal N}\left(0\right) $.

\end{theorem}

\begin{proof} Obviously $ {\mathcal M}\left(x\right)\subset{\mathcal N}\left(0\right) $. We have to show that $ {\mathcal M}\left(x\right)={\mathcal N}\left(0\right) $.
Let $ m\in{\mathcal N}\left(0\right) $. There exists an open $ {\mathcal V}\subset Y_{x} $, $ 0\in{\mathcal V} $, $ \varphi\in{\mathcal O}\left({\mathcal V}\right)\otimes M_{x} $ such that $ \partial\varphi=0 $ and
$ \varphi\left(0\right)=m $. We have to extend $ \varphi $ to some open set $ {\mathcal U}\subset X $. Let $ {\mathfrak g}=C_{{\mathfrak g}}\left(x\right)\oplus{\mathfrak l} $ as
$ {\mathfrak g}_{x} $-module. Define the map
\begin{equation}
\tau\colon {\mathfrak l}_{0} \times Y_{x} \to X
\notag\end{equation}
by the formula
\begin{equation}
\tau\left(l,y\right)=\exp  \operatorname{ad}\left(l\right) \left(x+y\right),
\notag\end{equation}
for any $ y\in Y_{x} $, $ l\in{\mathfrak l}_{0} $. Then $ \tau $ is a local isomorphism. Hence in some
neighborhood $ {\mathcal U}\subset X $, $ x\in{\mathcal U} $, $ x=\tau\left(l,y\right) $ and one can define
\begin{equation}
\psi\left(\tau\left(l,y\right)\right)=\exp  l \varphi\left(y\right).
\notag\end{equation}
Then $ \partial\psi=0 $ and $ \psi\left(x\right)=m $. Theorem is proven.\end{proof}

\section{Application to supercharacters }

The properties of $ M_{x} $ allow one to say something about the
superdimension and supercharacter of $ M $. First, we recall that $ \operatorname{sdim} M_{x}=\operatorname{sdim} M $.
Therefore

\begin{lemma} \label{lm106}\myLabel{lm106}\relax  If $ X_{M}\not=X $, then $ \operatorname{sdim} M=0 $. In particular, if a
finite-dimensional module $ M $ admits a central character whose degree of
atypicality is less than the defect of $ {\mathfrak g} $, then $ \operatorname{sdim} M=0 $.

\end{lemma}

Now let $ M $ be a finite-dimensional $ {\mathfrak g} $-module and $ h\in{\mathfrak h} $. Write
\begin{equation}
\operatorname{ch}_{M}\left(h\right)=\text{str}_{M}\left(e^{h}\right).
\notag\end{equation}
Obviously, $ \operatorname{ch}_{M} $ is $ W $-invariant analytic function on $ {\mathfrak h} $. We can write Taylor
series for $ \operatorname{ch}_{M} $ at $ h=0 $
\begin{equation}
\operatorname{ch}_{M}\left(h\right)=\sum_{i=0}^{\infty}p_{i}\left(h\right),
\notag\end{equation}
where $ p_{i}\left(h\right) $ is a homogeneous polynomial of degree $ i $ on $ {\mathfrak h} $. The order of $ \operatorname{ch}_{M} $
at zero is by definition the minimal $ i $ such that $ p_{i}\not\equiv 0 $.

\begin{theorem} \label{th9}\myLabel{th9}\relax  Assume that all odd roots of $ {\mathfrak g} $ are isotropic. Let $ M $ be a
finite-dimensional $ {\mathfrak g} $-module, $ s $ be the
codimension of $ X_{M} $ in $ X $. The order of $ \operatorname{ch}_{M} $ at zero is greater or equal than $ s $.
Moreover, the polynomial $ p_{s}\left(h\right) $ in Taylor series for $ \operatorname{ch}_{M} $ is determined
uniquely up to proportionality.

\end{theorem}

\begin{proof} The proof is based on the following Lemma, the proof of this
Lemma is similar to the proof of Lemma~\ref{lm2} (6). We leave it to the reader.

\begin{lemma} \label{lm153}\myLabel{lm153}\relax  Let $ x\in X $, $ h\in{\mathfrak g}_{0} $ and $ \left[h,x\right]=0 $. Then $ \operatorname{Ker} x $ and $ xM $ are
$ h $-invariant and str$ _{M}h= $str$ _{M_{x}}h $.

\end{lemma}

If $ X_{M}=X $, the statement of theorem is trivial. Let $ X_{M}\not=X $. By Theorem~\ref{th2}
there exists $ k $ less than the defect of $ {\mathfrak g} $ such that
\begin{equation}
X_{M}\subset\cup_{{\text A}\in S\text{, }|{\text A}|\leq k}\Phi\left({\text A}\right).
\notag\end{equation}
Let $ {\text A}=\left\{\alpha_{1},\dots ,\alpha_{k+1}\right\}\in S $, $ x=x_{1}+\dots +x_{k+1} $ for some nonzero $ x_{i}\in{\mathfrak g}_{\alpha_{i}} $. Then $ M_{x}=\left\{0\right\} $.
If $ h\in{\mathfrak h} $ satisfies $ \alpha_{1}\left(h\right)=\dots =\alpha_{k+1}\left(h\right)=0 $, then $ \left[h,x\right]=0 $. Hence by Lemma~\ref{lm153}
str$ _{M}h= $str$ _{M_{x}}h=0 $. Hence we just have proved the following property
\begin{equation}
\operatorname{ch}_{M}\left({\mathfrak h}_{{\text A}}^{\perp}\right)=0\text{ for all }{\text A}\in S,|{\text A}|=k+1.
\label{equ29}\end{equation}\myLabel{equ29,}\relax 
Let $ p_{i} $ be the
first non-zero polynomial in the Taylor series for $ \operatorname{ch}_{M} $ at zero. Then $ p_{i} $ also
satisfies~\eqref{equ29}. Let $ {\text B}=\left\{\alpha_{1},\dots ,\alpha_{k}\right\}\in S $ and $ \bar{p}_{i} $ be
the restriction of $ p_{i} $ to $ {\mathfrak h}_{{\text B}}^{\perp} $. If $ \bar{p}_{i}\not=0 $, then degree of $ \bar{p}_{i} $ is $ i $.
Since $ p_{i}\left({\mathfrak h}_{{\text B}\cup\alpha}^{\perp}\right)=0 $ for any $ \alpha\not=\pm\alpha_{i} $, $ \alpha\in{\text B}^{\perp} $, then $ \alpha $ divides $ \bar{p}_{i} $.
That gives the estimate on $ i $. Indeed, $ i $ is not less the number of all
possible $ \alpha $ up to proportionality, i.e. $ \frac{|\Delta_{1}\cap{\text B}^{\perp}|}{2}-|{\text B}| $. By Corollary~\ref{cor72}
the latter number is the codimension $ s $ of $ X_{M} $ in $ X $. Hence $ i\geq s $.

To prove the second statement we need to show
that if two homogeneous $ W $-invariant polynomials $ p $ and $ q $ of degree $ s $ satisfy~%
\eqref{equ29}, then $ p=cq $ for some $ c\in{\mathbb C} $. After restriction on $ {\mathfrak h}_{{\text B}}^{\perp} $
\begin{equation}
\bar{p}=a\Pi_{\alpha\in\left(\Delta^{+}\cap{\text B}^{\perp}\right)\backslash\pm{\text B}} \alpha\text{, }\bar{q}=b\Pi_{\alpha\in\left(\Delta^{+}\cap{\text B}^{\perp}\right)\backslash\pm{\text B}} \alpha
\notag\end{equation}
for some constants $ a $ and $ b $. Therefore there exists $ f=p-cq $ such that
$ f\left({\mathfrak h}_{{\text B}}^{\perp}\right)=0 $. Thus, $ f $ satisfies~\eqref{equ29} for $ k $ instead of $ k+1 $. Then the degree of $ f $
is bigger than $ s $, which implies $ f=0 $.\end{proof}

\section{Translation functor }

In this section we introduce translation functors, we use these
functors in the proof of Theorem~\ref{th3}. A translation functor is
a superanalogue of similar functor in category $ {\mathcal O} $ (see\cite{BGG76Cat}). For
superalgebras translation functors were used in \cite{Ser96KLpol} and \cite{Bru03KL}.

Let $ V $ be a $ {\mathfrak g} $-module, on which the center $ Z $ of the universal
enveloping algebra acts locally finitely. Then $ V=\oplus V^{\chi} $, where
\begin{equation}
V^{\chi}=\left\{v\in V \mid \left(z-\chi\left(z\right)\right)^{N}v=0, z\in Z\right\}.
\notag\end{equation}
Let $ {\mathcal B} $ be the category of all finitely generated $ {\mathfrak g} $-modules with finite
$ Z $-action. Then $ {\mathcal B} $ has a decomposition
\begin{equation}
{\mathcal B}=\oplus{\mathcal B}^{\chi},
\notag\end{equation}
where $ {\mathcal B}^{\chi} $ denotes the subcategory of all $ V\in{\mathcal B} $ such that $ V^{\chi}=V $.

Let $ E $ be a finite-dimensional $ {\mathfrak g} $-module. A {\em translation functor\/} $ T_{E}^{\chi} $ is
a functor in the category $ {\mathcal B} $, defined by
\begin{equation}
T_{E}^{\chi}\left(V\right)=\left(V\otimes E\right)^{\chi}.
\notag\end{equation}
To simplify the notation we also will write $ T_{E}^{\lambda} $ instead of $ T_{E}^{\chi_{\lambda}} $.

\begin{lemma} \label{lm7}\myLabel{lm7}\relax  $ T_{E}^{\chi} $ is an exact functor.

\end{lemma}

\begin{proof} Both tensoring with finite-dimensional vector space and the
projection on the component with a given central character are obviously
exact functors.\end{proof}

Denote by $ P\left(E\right) $ the set of all weights of $ E $ counted with
multiplicities.

\begin{lemma} \label{lm8}\myLabel{lm8}\relax 
\begin{enumerate}
\item
For the Verma module $ M_{\lambda} $, $ M_{\lambda}\otimes E $ has a finite
filtration $ \left\{0\right\}={\mathcal F}_{0}\subset\dots \subset{\mathcal F}_{q}=\left(M_{\lambda}\otimes E\right) $ of
length $ q=\dim  E $ such that $ {\mathcal F}_{i+1}/{\mathcal F}_{i} $ is ismorphic to $ M_{\lambda+\nu} $, $ \nu\in P\left(E\right) $;
\item
If $ V $ is a module generated by a highest vector of weight $ \lambda $, then
$ T_{E}^{\chi}\left(V\right) $ has a finite filtration $ \left\{0\right\}={\mathcal V}_{0}\subset\dots \subset{\mathcal V}_{r}=T_{E}^{\chi}\left(V\right) $ such that $ {\mathcal V}_{i}/{\mathcal V}_{i+1} $ is a
highest weight module of weight $ \lambda+\nu\in{\mathfrak h}_{\chi} $ for some $ \nu\in P\left(E\right) $.
\end{enumerate}
\end{lemma}

\begin{proof} The first statement can be found in \cite{BGG76Cat}. The second
one follows from the first and Lemma~\ref{lm7}.

{}\end{proof}

Let $ {\mathfrak b} $ be a Borel subalgebra of $ {\mathfrak g} $, $ V $ be a $ {\mathfrak g} $-module. A vector $ v\in V $ is
$ {\mathfrak b} $-{\em primitive\/} if $ {\mathfrak b}v\in{\mathbb C}v $.

\begin{lemma} \label{lm9}\myLabel{lm9}\relax  If $ v $ is a $ {\mathfrak b} $-primitive vector of $ \left(L_{\lambda}\otimes E\right) $ then the weight of
$ v $ equals $ \lambda+\nu $ for some $ \nu\in P\left(E\right) $.

\end{lemma}

\begin{proof} Introduce the order on $ {\mathfrak h}^{*} $ by putting $ \mu\leq\nu $ if $ \nu=\mu+\Sigma n_{\alpha}\alpha $
for some $ \alpha\in\Delta^{+} $ and $ n_{\alpha}\in{\mathbb Z}_{\geq0} $. Choose a maximal weight $ \gamma $ of $ L_{\lambda} $ such that
\begin{equation}
v=v_{1}\otimes w_{1}+\dots +v_{r}\otimes w_{r}+v'_{1}\otimes w'_{1}+\dots +v'_{t}\otimes w'_{t}
\notag\end{equation}
for some linearly independent $ v_{1},\dots ,v_{r}\in L_{\lambda} $ of weight $ \gamma $, linearly
independent weight vectors $ w_{1},\dots ,w_{r}\in E $ and some linearly independent weight
vectors $ v'_{1},\dots ,v'_{t}\in L_{\lambda} $ of weights different from $ \gamma $, $ w'_{1},\dots ,w'_{t}\in E $. For any
simple root element $ e\in{\mathfrak n} $ the condition $ e v=0 $ implies
\begin{equation}
e v_{1}\otimes w_{1}+\dots +e v_{r}\otimes w_{r}=0.
\notag\end{equation}
Since $ w_{1},\dots ,w_{r} $ are linearly independent, we must have $ e v_{i}=0 $. Therefore all
$ v_{i} $ are $ {\mathfrak b} $-primitive. But $ L_{\lambda} $ has a unique up to proportionality $ {\mathfrak b} $-primitive
vector. Therefore $ \gamma=\lambda $, $ r=1 $ and the weight of $ v $ is the sum of $ \lambda $ and the
weight of $ w_{1} $.\end{proof}

For any $ \lambda\in{\mathfrak h}^{*} $ put
\begin{equation}
{\mathfrak h}_{\lambda}={\mathfrak h}_{\chi_{\lambda}}\text{, }\Sigma_{\lambda}={\mathfrak h}_{\lambda}\cap\Sigma.
\notag\end{equation}

\begin{lemma} \label{lm10}\myLabel{lm10}\relax  Let $ \lambda,\mu\in\Sigma $ satisfy the conditions
\begin{equation}
\left(\lambda+P\left(E\right)\right)\cap\Sigma_{\mu}=\left\{\mu\right\}
\label{equ41}\end{equation}\myLabel{equ41,}\relax 
\begin{equation}
\left(\mu-P\left(E\right)\right)\cap\Sigma_{\lambda}=\left\{\lambda\right\}
\label{equ42}\end{equation}\myLabel{equ42,}\relax 
and $ \lambda $ is minimal in $ \left(\mu-P\left(E\right)\right)\cap{\mathfrak h}_{\lambda} $. Then
\begin{equation}
T_{E}^{\mu}\left(L_{\lambda}\right)=L_{\mu}\text{, }T_{E^{*}}^{\lambda}\left(L_{\mu}\right)=L_{\lambda}.
\notag\end{equation}
\end{lemma}

\begin{proof} By Lemma~\ref{lm8} $ \left(2\right) $ and~\eqref{equ41} $ T_{E}^{\mu}\left(L_{\lambda}\right) $ is a highest weight
module with highest weight $ \mu $. By Lemma~\ref{lm9} and~\eqref{equ41} $ T_{E}^{\mu}\left(L_{\lambda}\right) $ has a
unique up to proportionality $ {\mathfrak b} $-primitive vector. Therefore either $ T_{E}^{\mu}\left(L_{\lambda}\right)=L_{\mu} $
or $ T_{E}^{\mu}\left(L_{\lambda}\right)=\left\{0\right\} $. In the same way either $ T_{E}^{\lambda}\left(L_{\mu}\right)=L_{\lambda} $ or $ T_{E}^{\lambda}\left(L_{\mu}\right)=\left\{0\right\} $.

Our next observation is
\begin{equation}
\operatorname{Hom}_{{\mathfrak g}}\left(M\otimes E^{*},N\right)\cong\operatorname{Hom}_{{\mathfrak g}}\left(M,N\otimes E\right),
\label{equ94}\end{equation}\myLabel{equ94,}\relax 
hence, in particular
\begin{equation}
\operatorname{Hom}_{{\mathfrak g}}\left(T_{E^{*}}^{\lambda}\left(L_{\mu}\right),L_{\lambda}\right)\cong\operatorname{Hom}_{{\mathfrak g}}\left(L_{\mu},T_{E}^{\mu}\left(L_{\lambda}\right)\right).
\label{equ95}\end{equation}\myLabel{equ95,}\relax 
Therefore $ T_{E}^{\mu}\left(L_{\lambda}\right)=\left\{0\right\} $ iff $ T_{E}^{\mu}\left(L_{\mu}\right)=\left\{0\right\} $. Let us prove that $ T_{E}^{\mu}\left(L_{\lambda}\right)\not=\left\{0\right\} $. Note
that by Lemma~\ref{lm8} $ \left(1\right) $ and~\eqref{equ42}, $ T_{E^{*}}^{\lambda}\left(M_{\mu}\right) $ has a subquotient isomorphic to
$ M_{\lambda} $. Moreover, since $ \lambda $ is a minimal weight in $ \left(\mu-P\left(E\right)\right)\cap{\mathfrak h}_{\lambda} $, there is a
quotient in $ T_{E^{*}}^{\lambda}\left(M_{\mu}\right) $ isomorphic to $ M_{\lambda} $, hence there is a quotient isomorphic
to $ L_{\lambda} $. Therefore
\begin{equation}
\operatorname{Hom}_{{\mathfrak g}}\left(T_{E^{*}}^{\lambda}\left(M_{\mu}\right),L_{\lambda}\right)\not=\left\{0\right\}.
\notag\end{equation}
But then using~\eqref{equ94}
\begin{equation}
\operatorname{Hom}_{{\mathfrak g}}\left(M_{\mu},T_{E}^{\mu}\left(L_{\lambda}\right)\right)\cong\operatorname{Hom}_{{\mathfrak g}}\left(T_{E^{*}}^{\lambda}\left(M_{\mu}\right),L_{\lambda}\right)\not=0.
\notag\end{equation}
Therefore $ T_{E}^{\mu}\left(L_{\lambda}\right)\not=\left\{0\right\} $. Finally by~\eqref{equ95} $ T_{E^{*}}^{\lambda}\left(L_{\mu}\right)\not=\left\{0\right\} $.\end{proof}

\begin{lemma} \label{lm13}\myLabel{lm13}\relax  Let $ M $ be a finite-dimensional $ {\mathfrak g} $-module and $ N=T_{E}^{\chi}\left(M\right) $. Then
$ X_{N}\subset X_{M} $.

\end{lemma}

\begin{proof} Let $ x\in X\backslash X_{M} $, then $ M $ is free over $ {\mathbb C}\left[x\right] $, and $ M\otimes E $ is also free
over $ {\mathbb C}\left[x\right] $. Since $ N $ is a direct summand of $ M\otimes E $, then $ N $ is free over $ {\mathbb C}\left[x\right] $.
That implies $ x\notin X_{N} $.\end{proof}

\section{Reduction to the stable case }

Fix a set of simple roots and the Borel subalgebra $ {\mathfrak b}\subset{\mathfrak g} $ generated by $ {\mathfrak h} $
and simple roots. We say that a subalgebra $ {\mathfrak q} $ is {\em admissible\/} if $ {\mathfrak q} $ is generated
by $ {\mathfrak h} $, some subset of simple roots and their negatives. By $ \Delta\left({\mathfrak q}\right) $ we denote
the root system of $ {\mathfrak q} $. We call $ \lambda $ {\em stable\/} with respect to $ {\mathfrak q} $ if the following
conditions hold for any isotropic $ \alpha\in\Delta $, $ \left(\lambda+\rho,\alpha\right)=0 $ implies $ \alpha\in\Delta\left({\mathfrak q}\right) $.

In this section we assume that $ {\mathfrak g}={\mathfrak g}{\mathfrak l}\left(m|n\right) $. Then
\begin{equation}
\Delta_{0}=\left\{\varepsilon_{i}-\varepsilon_{j} \mid i,j\leq m\right\}\cup\left\{\delta_{i}-\delta_{j} \mid i,j\leq n\right\}\text{, }\Delta_{1}=\left\{\pm\left(\varepsilon_{i}-\delta_{j}\right) \mid i\leq m,j\leq n\right\}.
\notag\end{equation}
All odd roots are isotropic. The choice of the form on $ {\mathfrak h}^{*} $ is such that
$ \left(\varepsilon_{i},\varepsilon_{j}\right)=\delta_{ij} $, $ \left(\delta_{i},\delta_{j}\right)=-\delta_{ij} $. The defect $ d=\min \left(m,n\right) $. We choose a Borel
subalgebra $ {\mathfrak b} $ so that the simple roots are
\begin{equation}
\varepsilon_{1}-\varepsilon_{2},\dots ,\varepsilon_{m-1}-\varepsilon_{m},\varepsilon_{m}-\delta_{1},\dots ,\delta_{n-1}-\delta_{n}.
\notag\end{equation}
If $ \lambda+\rho=a_{1}\varepsilon_{1}+\dots +a_{m}\varepsilon_{m}+b_{1}\delta_{1}+\dots b_{n}\delta_{n} $, then $ \lambda\in\Sigma $ iff $ a_{i}-a_{i+1} $, $ b_{j}-b_{j+1}\in{\mathbb Z}_{>0} $ for
all $ i<m,j<n $. In other words, $ \lambda\in\Sigma $ iff $ \left<\lambda+\rho,\gamma^{\vee}\right>\in{\mathbb Z}_{>0} $ for all $ \gamma\in\Delta_{0}^{+} $. Since we
consider only atypical $ \lambda $ we may assume that $ a_{i},b_{j}\in{\mathbb Z} $.

\begin{lemma} \label{lm11}\myLabel{lm11}\relax  Let $ \lambda+\rho=a_{1}\varepsilon_{1}+\dots +a_{m}\varepsilon_{m}+b_{1}\delta_{1}+\dots b_{n}\delta_{n}\in\Sigma $.
If $ a_{i}+b_{j}=0 $ implies $ i>m-k $, then $ \lambda $ is stable for
$ {\mathfrak q} $ with simple roots $ \varepsilon_{m-k+1}-\varepsilon_{m-k+2},\dots ,\varepsilon_{m-1}-\varepsilon_{m},\varepsilon_{m}-\delta_{1},\dots ,\delta_{n-1}-\delta_{n} $.

\end{lemma}

\begin{proof} Trivial.\end{proof}

\begin{theorem} \label{th4}\myLabel{th4}\relax  Let $ {\mathfrak g}={\mathfrak g}{\mathfrak l}\left(m|n\right) $. If $ \lambda\in\Sigma $ and has the degree of
atypicality $ k $, then there exists a subalgebra $ {\mathfrak q}\subset{\mathfrak g} $ of defect $ k $,
a stable $ \mu\in\Sigma $ of degree atypicality $ k $ and translation functors
$ T_{1},\dots ,T_{r} $, $ T_{1}^{*},\dots ,T_{r}^{*} $ such that
\begin{equation}
L_{\mu}=T_{1}\dots T_{r}\left(L_{\lambda}\right)\text{, }L_{\lambda}=T_{r}^{*}\dots T_{1}^{*}\left(L_{\mu}\right).
\notag\end{equation}
\end{theorem}

\begin{proof} Translation functors which we use are always related with $ E $
being the standard representation or its dual. We will provide a
combinatorial algorithm, which constructs from a weight $ \lambda\in\Sigma $ a new weight
$ \mu\in\Sigma_{\lambda} $ in such way that $ \lambda $ and $ \mu $ satisfy the conditions of Lemma~\ref{lm10} and
therefore $ T_{E}^{\mu}\left(L\left(\lambda\right)\right)=L\left(\mu\right) $, $ T_{E^{*}}^{\lambda}\left(L\left(\mu\right)\right)=L\left(\lambda\right) $. Applying this algorithm several
times we obtain a sequence of weights $ \mu_{1},\dots ,\mu_{r} $ such that $ \mu_{r} $ is stable.
Let $ \lambda+\rho=a_{1}\varepsilon_{1}+\dots +a_{m}\varepsilon_{m}+b_{1}\delta_{1}+\dots b_{n}\delta_{n} $. Let $ g $ be maximal such
that $ a_{i}+b_{j}\not=0 $ for any $ i\leq g $, $ j\leq n $. If $ g=m-k $, then $ \lambda $ is stable as in Lemma~\ref{lm11}
and we can stop to apply the algorithm. Otherwise choose first $ i>g $ such that
$ a_{i}+b_{j}\not=0 $ for all $ j\leq n $. Construct $ \mu $ depending on the following
\begin{enumerate}
\item
If $ b_{j}\not=-a_{i}-1 $ for any $ j\leq n $, then put $ \mu=\lambda+\varepsilon_{i} $;
\item
If $ b_{j}=-a_{i}-1 $ for some $ j $ look at $ a_{i-1} $. If $ a_{i-1}=a_{i}+1 $, put $ \mu=\lambda+\delta_{j} $.
Otherwise go to the next step;
\item
If $ b_{j}=-a_{i}-1 $, $ a_{i-1}\not=a_{i}+1 $ find the maximal $ p $ such that $ b_{j+p}=b_{j}-p $. If
$ a_{i-1}+b_{j+p}>0 $, put $ \mu=\lambda-\delta_{j+p} $. Otherwise go to the next step.
\item
If $ a_{i-1}+b_{j+p}\leq0 $, then there exists $ t\leq p $ such that $ a_{i-1}+b_{j+t}=0 $. Put
$ \mu=\lambda-\varepsilon_{i-1} $.
\end{enumerate}

Note that at some point one always arrives to the case 2, that
decreases $ i $ and eventually increases $ g $. In the end one will come to the
stable weight.\end{proof}

Theorem~\ref{th4} and Lemma~\ref{lm13} imply

\begin{theorem} \label{th5}\myLabel{th5}\relax  Let $ {\mathfrak g}={\mathfrak g}{\mathfrak l}\left(m|n\right) $.
For any $ \lambda\in\Sigma $ there exists a stable $ \mu\in\Sigma $ with the same
degree of atypicality such that $ X_{L_{\lambda}}=X_{L_{\mu}} $.

\end{theorem}

\section{Proof of Theorem\protect ~\protect \ref{th3} for $ {\mathfrak g}{\mathfrak l}\left(m|n\right) $ }

In this section $ {\mathfrak g}={\mathfrak g}{\mathfrak l}\left(m|n\right) $ and $ \lambda $ is an integral dominant weight with degree
of atypicality $ k $. As Theorem~\ref{th2} is already proven we have to show only
that if $ {\text A}\in S $, $ |{\text A}|=k $, then $ \left(L_{\lambda}\right)_{x}\not=\left\{0\right\} $ for any $ x\in\Phi\left({\text A}\right) $. As follows from
Theorem~\ref{th5}, we may assume that $ \lambda $
is stable with respect to $ {\mathfrak q}={\mathfrak g}{\mathfrak l}\left(k|n\right) $. It is easy to check that $ \Phi\left({\text A}\right)\cap{\mathfrak q}\not=\varnothing $, and
therefore one may assume that $ x\in{\mathfrak q} $. On the other hand, $ L_{\lambda}=L_{\lambda}\left({\mathfrak q}\right)\oplus N $ as a module
over $ {\mathfrak q} $. Thus, it is sufficient to prove that $ \left(L_{\lambda}\left({\mathfrak q}\right)\right)_{x}\not=\left\{0\right\} $. In other words,
we reduce the theorem to the case of $ {\mathfrak g}{\mathfrak l}\left(k|n\right) $. Using the isomorphism
$ {\mathfrak g}{\mathfrak l}\left(k|n\right)\cong{\mathfrak g}{\mathfrak l}\left(n|k\right) $ we can repeat the above argument and reduce the theorem to
the case $ {\mathfrak g}={\mathfrak g}{\mathfrak l}\left(k|k\right) $. Summing up, Theorem~\ref{th3} is equivalent to the
following Lemma.

\begin{lemma} \label{lm100}\myLabel{lm100}\relax  Let $ {\mathfrak g}={\mathfrak g}{\mathfrak l}\left(k|k\right) $ and $ \lambda $ be an integral dominant weight with
degree of atypicality $ k $. Then $ \left(L_{\lambda}\right)_{x}\not=\left\{0\right\} $ for any $ x\in X $.

\end{lemma}

We prove Lemma~\ref{lm100} in several steps. We use the fact that $ {\mathfrak g} $ has the
$ {\mathbb Z} $-grading $ {\mathfrak g}={\mathfrak g}\left(-1\right)\oplus{\mathfrak g}\left(0\right)\oplus{\mathfrak g}\left(1\right) $ and $ {\mathfrak g}\left(-1\right),{\mathfrak g}\left(1\right) $ are irreducible components of $ X $. We
have $ k+1 $ open orbits on $ X $. Choose a representative $ x $ on each orbit in the
following way:
\begin{equation}
\left(
\begin{matrix}
0 & x^{+}
\\
x^{-} & 0
\end{matrix}
\right),
\notag\end{equation}
where $ x^{+} $ is the block matrix
\begin{equation}
\left(
\begin{matrix}
1_{p} & 0
\\
0 & 0
\end{matrix}
\right)
\notag\end{equation}
and $ x^{-} $ is the block matrix
\begin{equation}
\left(
\begin{matrix}
0 & 0
\\
0 & 1_{q}
\end{matrix}
\right);
\notag\end{equation}
here $ p+q=k $.

If $ x\in{\mathfrak g}\left(1\right) $, then $ x^{-}=0 $, if $ x\in{\mathfrak g}\left(-1\right) $, then $ x^{+}=0 $. In both cases the
stabilizer $ K $ of $ x $ in $ G_{0} $ is isomorphic to $ \operatorname{GL}\left(k\right) $ embedded diagonally in
$ G_{0}=\operatorname{GL}\left(k\right)\times\operatorname{GL}\left(k\right) $. By $ {\mathfrak k} $ we denote the Lie algebra of $ K $.

\begin{lemma} \label{lm101}\myLabel{lm101}\relax  If $ x\in{\mathfrak g}\left(\pm1\right) $ and $ M $ is a finite-dimensional $ {\mathfrak g} $-module, then $ M_{x} $
is a trivial $ K $-module.

\end{lemma}

\begin{proof} Follows from the fact $ C_{{\mathfrak g}}\left(x\right)=\left[x,{\mathfrak g}\right] $ and Lemma~\ref{lm204}.

{}\end{proof}

\begin{lemma} \label{lm151}\myLabel{lm151}\relax  Let $ {\mathfrak g}={\mathfrak g}{\mathfrak l}\left(m|n\right) $, $ {\mathfrak b} $ is the Borel subalgebra containing $ {\mathfrak g}\left(1\right) $,
$ M_{\lambda}=U\left({\mathfrak g}\right)\otimes_{U\left({\mathfrak b}\right)}C_{\lambda} $ be the Verma module. If $ \alpha $ is a negative isotropic root such
that $ \left(\lambda+\rho,\alpha\right)=0 $, then $ M_{\lambda} $ contains a $ {\mathfrak b} $-primitive vector of weight $ \lambda+\alpha $.

\end{lemma}

\begin{proof} Let $ I_{\alpha} $ be the set of all weights $ \lambda\in{\mathfrak h}^{*} $ such that $ M_{\lambda} $ has a
$ {\mathfrak b} $-primitive vector of weight $ \lambda+\alpha $. Then $ I_{\alpha} $ is Zariski closed, see for example
\cite{Dix96Env}. Let
\begin{equation}
H_{\alpha}=\left\{\lambda\in{\mathfrak h}^{*} \mid \left(\lambda+\rho,\alpha\right)=0\right\}.
\notag\end{equation}
We want to show that $ H_{\alpha}\subset I_{\alpha} $. Consider
\begin{equation}
H'_{\alpha}=\left\{\lambda\in H_{\alpha} \mid \left(\lambda+\rho,\beta\right)\not=0, \beta\not=\pm\alpha, \beta\in\Delta{\mathfrak g}_{1}\right\}.
\notag\end{equation}
It suffices to show that $ H'_{\alpha}\subset I_{\alpha} $. Consider $ {\mathfrak b}'={\mathfrak b}_{0}+{\mathfrak g}\left(-1\right) $. If $ v $ is a highest
vector of $ M_{\lambda} $ and $ X_{\beta}\in{\mathfrak g}_{\beta} $, then $ w=\Pi_{\beta\in\Delta\left({\mathfrak g}\left(-1\right)\right)}X_{\beta}v $ is a $ {\mathfrak b}' $-primitive, and
$ u=\Pi_{\beta\in\Delta\left({\mathfrak g}\left(1\right)\right)\backslash\left\{\alpha\right\}}X_{\beta}w $ is a $ {\mathfrak b} $-primitive. Since the weight of $ u $ equals $ \lambda+\alpha $, we
obtain $ H'_{\alpha}\subset I_{\alpha} $ as required.\end{proof}

The $ {\mathbb Z} $-grading on $ {\mathfrak g} $ induces the $ {\mathbb Z} $-grading on an irreducible $ {\mathfrak g} $-module
$ M=M\left(0\right)\oplus M\left(-1\right)\oplus\dots \oplus M\left(-k^{2}\right) $ in the following way
\begin{equation}
M\left(0\right)=\operatorname{Ker} {\mathfrak g}\left(1\right)\text{, }M_{i}={\mathfrak g}\left(-1\right)M\left(i+1\right).
\notag\end{equation}
Each $ M\left(i\right) $ is a $ {\mathfrak g}_{0} $-submodule of $ M $.

\begin{lemma} \label{lm104}\myLabel{lm104}\relax  Let $ x\in{\mathfrak g}\left(\pm1\right) $, $ M=L_{\lambda} $ for a dominant integral $ \lambda $ of degree
atypicality $ k $. Then $ M\left(0\right) $ contains one trivial $ K $-submodule and $ M\left(-1\right) $
does not have trivial $ K_{x} $-submodules.

\end{lemma}

\begin{proof} Since the degree of atypicality of $ \lambda $ is $ k $, one can write
\begin{equation}
\lambda=a_{1}\varepsilon_{1}+\dots +a_{k}\varepsilon_{k}-a_{k}\delta_{1}-\dots -a_{1}\delta_{1}.
\notag\end{equation}
We denote by $ V\left(a_{1},\dots ,a_{k}\right) $ the irreducible $ {\mathfrak g}{\mathfrak l}\left(k\right) $-module with highest weight
$ \left(a_{1},\dots ,a_{k}\right) $ and by $ L_{\lambda}\left({\mathfrak g}_{0}\right) $ the irreducible $ {\mathfrak g}_{0} $-module with highest
weight $ \lambda $. Since $ M\left(0\right) $ is isomorphic to $ L_{\lambda}\left({\mathfrak g}_{0}\right) $, then
\begin{equation}
M\left(0\right)\cong V\left(a_{1},\dots ,a_{k}\right)\otimes V^{*}\left(a_{1},\dots ,a_{k}\right)
\notag\end{equation}
as $ K $-module, which has exactly one trivial component. Hence the first
statement is true.

Obviously $ M\left(-1\right) $ is a submodule in
\begin{equation}
L_{\lambda}\otimes{\mathfrak g}\left(-1\right)\subset\oplus_{\alpha\in\Delta\left({\mathfrak g}\left(-1\right)\right)}L_{\lambda+\alpha}\left({\mathfrak g}_{0}\right).
\notag\end{equation}
However, $ \left(\lambda+\rho, \varepsilon_{i}-\delta_{k+1-i}\right)=0 $, therefore by Lemma~\ref{lm151} $ M\left(-1\right) $ does not contain
the component $ L_{\lambda+\delta_{k+1-i}-\varepsilon_{i}}\left({\mathfrak g}_{0}\right) $ for all $ i=1,\dots ,k $. Hence $ M\left(-1\right) $ is a
$ K $-submodule of the $ K $-module
\begin{equation}
\oplus_{i\not=j}V\left(a_{1},\dots ,a_{i}-1,\dots ,a_{k}\right)\otimes V^{*}\left(a_{1},\dots ,a_{j}-1,\dots ,a_{k}\right).
\notag\end{equation}
Therefore $ M\left(-1\right) $ does not contain $ K $-trivial submodules.\end{proof}

\begin{lemma} \label{lm105}\myLabel{lm105}\relax  Let $ x\in{\mathfrak g}\left(\pm1\right) $ belong to an open $ G_{0} $-orbit, $ M=L_{\lambda} $ for a
dominant integral $ \lambda $ of degree of atypicality $ k $ and $ N $ be a trivial $ K $-submodule
in $ M\left(0\right) $. Then $ N\subset M_{x} $ and therefore $ M_{x}\not=\left\{0\right\} $.

\end{lemma}

\begin{proof} If $ x\in{\mathfrak g}\left(1\right) $, then $ xN=0 $. Since $ x:M\left(-1\right) \to M\left(0\right) $ is a
homomorphism of
$ {\mathfrak g}_{x} $-modules and $ M\left(-1\right) $ does not have trivial $ {\mathfrak g}_{x} $-submodules, then $ N $ does not
belong to $ \operatorname{Im} x $. If $ x\in{\mathfrak g}\left(-1\right) $, then $ N $ clearly is not in $ \operatorname{Im} x $. Since $ x:M\left(0\right) \to
M\left(-1\right) $ is a homomorphism of $ {\mathfrak g}_{x} $-modules and $ M\left(-1\right) $ does not contain trivial
$ {\mathfrak g}_{x} $-submodules, then $ xN=0 $.\end{proof}

Lemma~\ref{lm105} shows that $ M_{x}\not=\left\{0\right\} $ in two special cases: $ x\in{\mathfrak g}\left(1\right) $ or $ x\in{\mathfrak g}\left(-1\right) $.
Now we will show the same for each open $ G_{0} $-orbit on $ C $. Let $ y_{c} $ be an odd
element in $ {\mathfrak g} $ given by
\begin{equation}
\left(
\begin{matrix}
0 & y_{c}^{+}
\\
y_{c}^{-} & 0
\end{matrix}
\right),
\notag\end{equation}
where $ y_{c}^{+} $ is the block matrix
\begin{equation}
\left(
\begin{matrix}
1_{p} & 0
\\
0 & c1_{q}
\end{matrix}
\right)
\notag\end{equation}
and $ y_{c}^{-} $ is the block matrix
\begin{equation}
\left(
\begin{matrix}
c1_{p} & 0
\\
0 & 1_{q}
\end{matrix}
\right);
\notag\end{equation}
here $ p+q=k $, $ c\in{\mathbb C} $. Note that $ y_{c}\notin X $ if $ c\not=0 $, but $ \left[y_{c},y_{c}\right] $ lies in the center of $ {\mathfrak g} $.
If $ M=L_{\lambda} $ has the degree atypicality $ k $, then the center of $ {\mathfrak g} $ acts by zero on
$ M $. Hence $ M_{y_{c}}=\operatorname{Ker} y_{c}/\operatorname{Im} y_{c} $ is well defined. Lemma~\ref{lm105} implies $ M_{y_{1}}\not=\left\{0\right\} $. If
$ c\not=0 $, then there exists $ g\in G_{0} $ such that $ y_{c}=c^{1/2} \operatorname{Ad}_{g}\left(y_{1}\right) $. Therefore
$ M_{y_{c}}\not=\left\{0\right\} $ for any $ c\not=0 $. The continuity argument shows that $ M_{y_{0}}\not=\left\{0\right\} $. But $ y_{0}\in X $ is
an element on an open orbit. Therefore Lemma~\ref{lm100} and Theorem~\ref{th3} are
proven.

\section{Application to $ H\left({\mathfrak g}\left(-1\right);M\right) $ for $ {\mathfrak g}{\mathfrak l}\left(m|n\right) $ }

Let $ {\mathfrak g}={\mathfrak g}{\mathfrak l}\left(m|n\right) $, then $ {\mathfrak g}\left(-1\right) $ is an abelian subalgebra and the cohomology
$ H\left({\mathfrak g}\left(-1\right);M\right) $ determine the character of a finite-dimensional module $ M $. On the
other hand, $ {\mathfrak g}\left(-1\right) $ is an irreducible component of $ X $. The complex calculating
$ H\left({\mathfrak g}\left(-1\right);M\right) $ is
\begin{equation}
\partial:{\mathcal O}\left({\mathfrak g}\left(-1\right)\right)\otimes M \to {\mathcal O}\left({\mathfrak g}\left(-1\right)\right)\otimes M,
\notag\end{equation}
where $ \partial $ is the same as for the sheaf $ {\mathcal M} $. One can consider the
localization of this complex and the corresponding coherent sheaf $ {\mathcal H}_{M} $ is
the restriction of $ {\mathcal M} $ on $ {\mathfrak g}\left(-1\right) $.

Theorem~\ref{th3} and Theorem~\ref{th211} imply the following

\begin{theorem} \label{th99}\myLabel{th99}\relax  Let $ M $ be an irreducible finite-dimensional module with
central character $ \chi $ and the degree of atypicality of $ \chi $ equal $ k $. Then $ \operatorname{supp}
{\mathcal H}_{M}=\bar{X}_{k}\cap{\mathfrak g}\left(-1\right) $.

\end{theorem}

\begin{lemma} \label{lm99}\myLabel{lm99}\relax  Let $ M $ be a typical finite-dimensional module. Then $ \operatorname{supp}
{\mathcal H}_{M}=\left\{0\right\} $ and $ {\mathcal H}_{M}\left(0\right)=H^{0}\left({\mathfrak g}\left(-1\right),M\right) $.

\end{lemma}

\begin{proof} Since $ M $ is typical, then $ M $ is a free $ {\mathfrak g}\left(-1\right) $ module and
$ H^{i}\left({\mathfrak g}\left(-1\right),M\right)=0 $ for $ i>0 $.\end{proof}

\begin{theorem} \label{th100}\myLabel{th100}\relax  Let $ x\in X_{k}\cap{\mathfrak g}\left(-1\right) $, $ M=L_{\lambda} $, the degree of
atypicality of $ \lambda $ be $ k $, and $ Z=G_{0}x $. Then $ {\mathcal H}_{M}\left(Z\right) $ is the sheaf of section of
the $ G_{0} $-vector bundle inuced by $ \left({\mathfrak g}_{x}\right)_{0} $-module $ H^{0}\left({\mathfrak g}_{x}\cap{\mathfrak g}\left(-1\right);M_{x}\right) $.

\end{theorem}

\begin{proof} Follows from Lemma~\ref{lm99}, Theorem~\ref{th211} and Theorem~\ref{th209}.
\end{proof}

\bibliography{ref,outref,mathsci}

\def\cprime{$'$} \def\cprime{$'$} \def\cprime{$'$} \def\cprime{$'$}
  \def\cprime{$'$}
\providecommand{\bysame}{\leavevmode\hbox to3em{\hrulefill}\thinspace}
\providecommand{\MR}{\relax\ifhmode\unskip\space\fi MR }
\providecommand{\MRhref}[2]{%
  \href{http://www.ams.org/mathscinet-getitem?mr=#1}{#2}
}
\providecommand{\href}[2]{#2}
\begin{thebibliography}{10}

\bibitem{Kac77Lie}
V.~G. Kac, \emph{Lie superalgebras}, Adv. Math. \textbf{26} (1977), 8--96.

\bibitem{Dix96Env}
J.~Dixmier, \emph{Enveloping algebras}, Graduate Studies in Mathematics,
  vol.~11, American Mathematical Society, Providence, RI, 1996, Revised reprint
  of the 1977 translation.

\bibitem{Kac78Rep}
V.~G. Kac, \emph{Representations of classical {L}ie superalgebras}, Lecture
  {N}otes in {M}ath. \textbf{676} (1978), 597--626.

\bibitem{Ser99Inv}
A.~Sergeev, \emph{The invariant polynomials on simple {L}ie superalgebras},
  Representation theory \textbf{3} (1999), 250--280.

\bibitem{Gru00Coh}
C.~Gruson, \emph{Sur la cohomologie des super alg\`ebres de {L}ie \'etranges},
  Transform. Groups \textbf{5} (2000), no.~1, 73--84.

\bibitem{BGG76Cat}
J.~Bernstein, I.~M. Gelfand, and S.I. Gelfand, \emph{Category of {G}--modules},
  Funct. Anal. and Appl. \textbf{10} (1976), 87--92.

\bibitem{Vog91Ass}
David~A. Vogan, Jr., \emph{Associated varieties and unipotent representations}.

\bibitem{Knop05inv}
Friedrich Knopp, \emph{Invariant functions on symplectic representations},
  preprint, arXiv:math.AG/0506171, Jun 2005.

\bibitem{Gor04Kac}
Maria Gorelik, \emph{The {K}ac construction of the centre of $u(g)$ for {L}ie
  superalgebras}, J. Nonlinear Math. Phys. (2004).

\bibitem{Bru03KL}
Jonathan Brundan, \emph{Kazhdan-lusztig polynomials and character formulae for
  the {L}}.

\bibitem{KatOchi01Deg}
Shohei Kato and Hiroyuki Ochiai, \emph{The degrees of orbits of the
  multiplicity-free actions}, Ast\'erisque (2001), no.~273, 139--158, Nilpotent
  orbits, associated cycles and Whittaker models for highest weight
  representations.

\bibitem{Kac84Lap}
Victor~G. Kac, \emph{Laplace operators of infinite-dimensional {L}ie algebras
  and theta functions}, Proceedings of the National Academy of Sciences of the
  United States of America \textbf{81} (1984), no.~2, 645--647.

\bibitem{NishOchiTan01Ber}
Kyo Nishiyama, Hiroyuki Ochiai, and Kenji Taniguchi, \emph{Bernstein degree and
  associated cycles of {H}arish-{C}handra modules---{H}ermitian symmetric
  case}, Ast\'erisque (2001), no.~273, 13--80, Nilpotent orbits, associated
  cycles and Whittaker models for highest weight representations.

\bibitem{Gru00Ide}
Caroline Gruson, \emph{Sur l'id\'eal du c\^one autocommutant des super
  alg\`ebres de {L}ie basiques classiques et \'etranges}, Ann. Inst. Fourier
  (Grenoble) \textbf{50} (2000), no.~3, 807--831.

\bibitem{Gru03Coh}
\bysame, \emph{Cohomologie des modules de dimension finie sur la super
  alg\`ebre de {L}ie {$\germ{osp}(3,2)$}}, J. Algebra \textbf{259} (2003),
  no.~2, 581--598.

\bibitem{KacWak94Aff}
V.~G. Kac and M.~Wakimoto, \emph{Integrable highest weight modules over affine
  superalgebras and number theory}, Progress in Math. \textbf{123} (1994),
  415--456.

\bibitem{Ser96KLpol}
V.~Serganova, \emph{Kazhdan-{L}usztig polynomials and character formula for the
  {L}ie superalgebra ${\germ g}{\germ l}(m\vert n)$}, Selecta Math. (N.S.)
  \textbf{2} (1996), no.~4, 607--651.

\end{thebibliography}
\end{document}